\documentclass[reqno,a4paper,12pt]{amsart} 
\usepackage{ifthen,latexsym,amssymb,amsmath,bbm,fixmath}

\usepackage{tikz}
\usepackage{hyperref}

\usepackage{thmtools}

\usepackage[shortlabels]{enumitem} 

\renewcommand{\ge}{\geqslant}
\renewcommand{\le}{\leqslant}

\newif\ifnotesw\noteswtrue

\newcommand{\hide}[1]{}

\newcommand{\beq}[1]{\begin{equation}\label{#1}}
\newcommand{\eeq}{\end{equation}}

\newtheorem{theorem}{Theorem}[section]

\newcommand{\bpf}[1][Proof.]{\smallskip\noindent{\it #1} }
\newcommand{\epf}{\qed \medskip}

\newtheorem{claim}{Claim}

\newtheorem{remark}[theorem]{Remark}

\newtheorem{proposition}[theorem]{Proposition}

\newtheorem{definition}[theorem]{Definition}

\newcommand{\comm}[1]{}

\usepackage[capitalize, nameinlink]{cleveref}
\crefname{theorem}{Theorem}{Theorems}
\crefname{proposition}{Proposition}{Propositions}
\crefname{observation}{Observation}{Observations}
\crefname{lemma}{Lemma}{Lemmas}
\crefname{claim}{Claim}{Claims}
\crefname{problem}{Problem}{Problems}
\crefname{conjecture}{Conjecture}{Conjectures}
\crefname{question}{Question}{Questions}
\crefname{example}{Example}{Examples}
\crefname{fact}{Fact}{Facts}

\renewcommand{\deg}{\operatorname{deg}}
\newcommand{\dom}{\operatorname{dom}}

\newcommand{\dist}{\operatorname{dist}}

\newcommand{\fB}{\mathcal{B}}
\newcommand{\fC}{\mathcal{C}}
\newcommand{\fD}{\mathcal{D}}
\newcommand{\fE}{\mathcal{E}}

\newcommand{\fG}{\mathcal{G}}
\newcommand{\fH}{\mathcal{H}}

\newcommand{\fP}{\mathcal{P}}

\newcommand{\fV}{\mathcal{V}}

\newcommand{\one}{{\bf 1}}

\newcommand{\Path}[4]{\ifthenelse{\equal{#1}{}}{P_{#2}(#3,#4)}{P_{#2}(#3,#4,#1)}}

\usepackage[shortlabels]{enumitem}

\begin{document}
\title{Measurable Vizing's theorem}

	\author{Jan Greb\'ik}
	\address{University of Warwick Mathematics Institute, Coventry CV4 7AL, UK.}
	\email{jan.grebik@warwick.ac.uk}

\maketitle

\begin{abstract}
    We prove a \emph{full} measurable version of Vizing's theorem for bounded degree Borel graphs, that is, we show that every Borel graph $\fG$ of degree uniformly bounded by $\Delta\in \mathbb{N}$ defined on a standard probability space $(X,\mu)$ admits a $\mu$-measurable proper edge coloring with $(\Delta+1)$-many colors.
    This answers a question of Marks \emph{[Question 4.9, J. Amer. Math. Soc. 29 (2016)]} also stated in Kechris and Marks as a part of \emph{[Problem 6.13, survey (2020)]}, and extends the result of the author and Pikhurko \emph{[Adv. Math. 374, (2020)]} who derived the same conclusion under the additional assumption that the measure $\mu$ is $\fG$-invariant. 
\end{abstract}

\section{Introduction}

\emph{Vizing's theorem} is a fundamental result in graph theory that relates the number of colors needed to properly color edges of a given graph $G$, the so-called \emph{chromatic index $\chi'(G)$ of $G$}, with its maximum degree; it states that if the maximum degree of $G$ is $\Delta\in \mathbb{N}$, then $\chi'(G)\le \Delta+1$.
Together with \emph{K\"{o}nig's line coloring theorem} (that states that $\chi'(G)=\Delta$ under the additional assumption that $G$ is bipartite), these classical results laid the foundation of edge-colouring, an important and active area of graph theory; see, for example, the recent book on edge-colouring by Stiebitz, Scheide, Toft and Favrholdt \cite{EdgeColoringBook}.

In this paper we study Vizing's theorem from the perspective of \emph{measurable graph combinatorics}, a subfield of descriptive set theory that lies at the intersection of measure theory, random processes, dynamics, group theory, combinatorics and distributed computing; a (non-exhaustive) sample of results related to the field (and to our investigation) includes \cite{laczk, marksunger,measurablesquare,doughertyforeman,marks2016baire,gaboriau,KST,DetMarks,millerreducibility,AsymptoticDim,csokagrabowski,grebik2020measurable,Bernshteyn2021LLL,brandt_chang_grebik_grunau_rozhon_vidnyaszkyhomomorphisms}.
In its most abstract form, measurable graph combinatorics systematically studies the existence of \emph{measurable} or \emph{definable} solutions to various graph coloring problems defined on so-called \emph{Borel graphs}, graphs where the vertex set $(V,\fB)$ is a \emph{standard Borel space}, for example the unit interval, and the edge relation $E$ is a Borel subset of $[V]^2$, where $[V]^2$ denotes the set of unordered pairs endowed with the canonical Borel structure coming from $V$.
This study originated in the seminal paper of Kechris, Solecki and Todor\v{c}evi\' c \cite{KST} and since then found many applications in various areas of central mathematics, see the surveys of Kechris and Marks \cite{kechris_marks2016descriptive_comb_survey} and Pikhurko \cite{pikhurko2021descriptive_comb_survey}.

Basic questions about vertex colorings of bounded degree Borel graphs are well-understood.
Recall that the \emph{Borel chromatic number $\chi_\fB(\fG)$} of a Borel graph $\fG=(V,\fB,E)$ is the smallest $k\in \mathbb{N}$ such that there is a decomposition $V=V_1\sqcup\dots\sqcup V_k$, where $V_i$ is a Borel subset of $V$ that does not span any edge of $\fG$ for every $i\in [k]$.
In the paper \cite{KST} Kechris, Solecki and Todor\v{c}evi\'{c} showed that there is a \emph{Borel measurable} version of the classical greedy algorithm for proper vertex coloring, thus proved that $\chi_\fB(\fG)\le \Delta+1$ for every Borel graph $\fG=(V,\fB,E)$ of maximum degree $\Delta\in \mathbb{N}$.
In the groundbreaking paper \cite{DetMarks}, Marks found an example of acylic $\Delta$-regular Borel graph $\fG$ such that $\chi_\fB(\fG)=\Delta+1$ for every $\Delta\in \mathbb{N}$, thus concluding that there is no Borel analogue of the \emph{classical Brooks' theorem} from finite combinatorics.
On the other hand, Conley, Marks and Tucker-Drob \cite{BrooksMeas} proved that Brooks' theorem holds once \emph{Borel measurability} is relaxed either in the sense of measure theory or topology.
For example, if $\fG=(V,\fB,E)$ is a Borel graph of degree bounded by $\Delta\ge 3$ that does not contain a clique on $\Delta+1$ vertices and $\mu$ is a \emph{Borel probability measure on $(V,\fB)$}, then there is a $\mu$-null Borel set $Y$ such that $\chi_{\fB}(\fG\upharpoonright (V\setminus Y))=\Delta$.
Recently these results were refined by Bernshteyn \cite{Bernshteyn2021LLL} and Brandt, Chang, the author, Grunau, Rozho\v n, and Vidny\' anszky \cite{brandt_chang_grebik_grunau_rozhon_vidnyaszkyhomomorphisms} using ideas from the theory of distributed computing.

In contrast, similar questions about \emph{edge colorings} are not yet fully understood.
A \emph{Borel matching} of a Borel graph $\fG=(V,\fB,E)$ is a Borel subset $M\subseteq E$ that is a matching, i.e., every two edges $e, f\in M$ are either equal (as unordered pairs) or do not share a vertex.
The \emph{Borel chromatic index $\chi'_B(\fG)$} is defined to be the smallest $k\in \mathbb{N}$ such that there is a decomposition $E=E_1\sqcup\dots\sqcup E_k$, where $E_i$ is a Borel matching for every $i\in [k]$.
Equivalently, this can be stated in the language of edge colorings as finding the smallest $k\in \mathbb{N}$ such that there is a Borel map $c:E\to [k]$ that is a \emph{proper edge coloring}, i.e., $c(e)\not=c(f)$ whenever $e\not =f\in E$ share a vertex.
The greedy upper bound \cite{KST} implies that $\chi'_B(\fG)\le 2\Delta-1$ for every Borel graph $\fG$ of maximum degree bounded by $\Delta\in \mathbb{N}$.
In the same paper \cite{DetMarks}, Marks found an example of acylic $\Delta$-regular Borel graph $\fG$ such that $\chi'_B(\fG)=2\Delta-1$ for every $\Delta\in \mathbb{N}$.
This shows that Borel analogue of Vizing's theorem does not hold.
Similarly as in the case of vertex colorings, it is natural to ask if Vizing's theorem hold once Borel measurability is relaxed either in the sense of measure theory or topology.
This question (for both measure and topology) is explicitly stated in the survey of Kechris and Marks \cite[Problem 6.13]{kechris_marks2016descriptive_comb_survey}.
Marks asked about the measurable relaxation for regular graphs \cite[Question 4.9]{DetMarks}, and, in fact, the same question for invariant measures, so-called \emph{graphings}, was raised earlier by Ab\' ert \cite{AbertQuestions}.

In this paper we completely resolve the measurable case, see \cref{thm:MainIntro}.
Before we state the result we formalize the definitions and discuss relevant results from recent years.
Recall that we always assume that the Borel graph $\fG=(V,\fB,E)$ in question has degree uniformly bounded by $\Delta\in \mathbb{N}$.
Given a \emph{Borel probability measure} $\mu$ on $(V,\fB)$, we define the \emph{$\mu$-chromatic index of $\fG$, $\chi_\mu'(\fG)$,} to be the minimal $k\in \mathbb{N}$ such that there is a Borel map $c:E\to [k]$ that satisfies
$$\mu\left(\left\{v\in V:c \ \operatorname{is \ not \ proper \ at \ } v\right\}\right)=0,$$
where $c$ \emph{is not proper at $v\in V$} if there are edges $e\not =f \in E$ that are adjacent to $v$ and $c(e)=c(f)$.
In another words, $c$ is a proper edge coloring at $\mu$-almost every vertex $v\in V$.
We say that $\mu$ is \emph{$\fG$-invariant}, if $\mu(g(C))=\mu(C)$ for every $C\in \fB$ and every Borel injection $g:C\to V$ that satisfies that $x$ and $g(x)$ are in the same connected component of $\fG$ for every $x\in C$.
In that case the quadruple $\fG=(V,\fB,\mu,E)$ is called a \emph{graphing}.

Cs\'{o}ka, Lippner and Pikhurko \cite{csokalippnerpikhurko} showed that $\chi'_\mu(\fG)\le \Delta+1$ for a graphing $\fG=(V,\fB,\mu,E)$ that does not contain odd cycles and proved an upper bound of $\Delta+O(\sqrt{\Delta})$ colors for graphings in general.
In a related result, Bernshteyn \cite[Theorem 1.3]{BernshteynEarlyLLL} proved that $\Delta+o(\Delta)$ colors are enough (even for the so-called list-colouring version) provided that the graphing factors to the shift action $\Gamma\curvearrowright [0, 1]^\Gamma$ of a finitely generated group $\Gamma$.
Answering the question of Ab\' ert, the author and Pikhurko \cite{grebik2020measurable} proved a measurable version of Vizing's theorem for graphings, that is, $\chi'_\mu(\fG)\le \Delta+1$ for any graphing $\fG=(V,\fB,\mu,E)$.
Interestingly, the technique developed in \cite{grebik2020measurable} was greatly extended by Bernshteyn \cite{BernshteynVizing} who found striking applications to the LOCAL model of distributed computing, see \cite{BernshteynVizing,Christiansen,BernshteynVizing2} for a current development in that direction.

Bernshteyn \cite{BernshteynEarlyLLL}, the author and Pikhurko \cite{grebik2020measurable}, T\' oth \cite{TothShcreier} and the author \cite{grebikApprox} investigated weaker notion of \emph{approximate} edge colorings.
In this setting, we require to find for each probability measure $\mu$ and $\epsilon>0$ a coloring of edges that is not correct or undefined for at most $\epsilon$ fraction of all edges.
Here, the analogues of Vizing's theorem as well as K\"{o}nig's line coloring theorem hold, see \cite{grebik2020measurable} and \cite{grebikApprox}.

Bowen and Weilacher \cite{BowenWeilacher} investigated Vizing's theorem in the context of (Borel) asymptotic dimension introduced in \cite{AsymptoticDim}.
As an application they derived that $\Delta+1$ colors are enough for edge colorings of bipartite graphs in the topological relaxation sense and for measures that are hyperfinite.
Stronger results in the special case of free Borel $\mathbb{Z}^d$-actions, that is, Borel edge coloring with $2d$ colors which is the analogue of K\"{o}nig's line coloring theorem in this setting, were obtained independently around the same time in \cite{ArankaFeriLaci,spencer_personal,grebik_rozhon2021toasts_and_tails,weilacher2021borel}.

Recently, Qian and Weilacher \cite{QianWeilacher} found connections of the topological relaxation to computable combinatorics which allowed them to derive an upper bound of $\Delta+2$ colors for the Baire measurable analogue of Vizing's theorem, the full topological analogue, that is $\Delta+1$ colors only, remains an interesting open problem.

Next, we state our main result, which is the full analogue of Vizing's theorem in the measurable setting.

\begin{theorem}\label{thm:MainIntro}
    Let $(V,\fB)$ be a standard Borel space, $\Delta\in \mathbb{N}$, $\fG=(V,\fB,E)$ be a Borel graph of uniformly bounded degree by $\Delta\in \mathbb{N}$ and $\mu$ be a Borel probability measure on $(V,\fB)$.
    Then $\chi_\mu'(\fG)\le \Delta+1$, i.e., there is a Borel map $c:E\to [\Delta+1]$ that is a proper edge coloring at $\mu$-almost every vertex $v\in V$.
\end{theorem}

In the proof we combine the technique of augmenting iterated Vizing chains introduced by the author and Pikhurko in \cite{grebik2020measurable} together with the following result, \cref{thm:CocycleIntro}, that is interesting in its own right and might be useful for other applications as well.\footnote{After the first version of this paper appeared on arXiv, we were informed by G\'{a}bor Elek that Gabriella Kuhn \cite[Lemma~1]{KuhnCocycle} proved the same result in the context of group actions, see \cref{rem:Kuhn}.}
Before we state the result, we need several definitions.
Let $\fG=(V,\fB,E)$ be a Borel graph and $\mu$ be a Borel probability measure on $(V,\fB)$.
Recall that $\mu$ is called $\fG$-quasi-invariant if $\mu([A]_\fG)=0$, whenever $\mu(A)=0$, where $[A]_\fG$ is the union of all connected component of $\fG$ that have non-empty intersection with $A$.
It is a standard fact, see eg \cite[Chapter~II, Section~8]{topics}, that if $\mu$ is $\fG$-quasi-invariant, then there is a function $\rho_\mu$, called the \emph{Radon--Nikodym cocycle of $\fG$}, that takes values in $(0,+\infty)$, is defined for every ordered pair of points $x,y\in V$ that are in the same connected component of $\fG$ and satisfies the following \emph{mass transport principle}:
$$\mu(g(C))=\int_C \rho_{\mu}(x,g(x)) \ d\mu$$
for every $C\subseteq V$ and an injective Borel map $g:C\to X$ that satisfies for every $x\in C$ that $g(x)$ is in the same connected component of $\fG$ as $x$.
While in general $\rho_\mu$ can behave chaotically, the next result shows that one can always pass to an equivalent measure $\nu$ such that $\rho_\nu$ is bounded on edges of $\fG$.

\begin{theorem}[see also \cite{KuhnCocycle} Lemma~1]\label{thm:CocycleIntro}
    Let $(V,\fB)$ be a standard Borel space, $\Delta\in \mathbb{N}$, $\fG=(V,\fB,E)$ be a Borel graph of uniformly bounded degree by $\Delta\in \mathbb{N}$ and $\mu$ be a Borel probability measure on $(V,\fB)$ that is $\fG$-quasi-invariant.
    Then there is an equivalent Borel probability measure $\nu$ on $(V,\fB)$ such that
    $$\frac{1}{4\Delta}\le \rho_\nu(x,y)\le 4\Delta$$
    for every edge $(x,y)\in E$.

    In particular, if $\dist_\fG(x,y)=k\in \mathbb{N}$, then $\rho_\nu(x,y)\le (4\Delta)^k$, where $\dist_\fG$ is the graph distance on $\fG$.
\end{theorem}

\begin{remark}\label{rem:Kuhn}
    Kuhn \cite[Lemma 1]{KuhnCocycle} showed that if a countable group $\Gamma$ acts in a quasi-invariant fashion on a standard probability space $(X,\mu)$, then there is an equivalent measure $\nu$ such that the cocycle $\rho_\nu({-},\gamma\cdot {-})$ (as a function $X\to (0,\infty)$) is bounded for every fixed $\gamma\in \Gamma$.
    This result combined with the fact that every Borel graph $\fG$ of degree bounded by $\Delta(\fG)<+\infty$ can be generated by $2\Delta(\fG)-1$ involutions (which follows from the Lusin--Novikov uniformization theorem \cite[Theorem~18.10]{kechrisclassical}) implies \cref{thm:CocycleIntro} possibly with a constant larger than $4\Delta$.
\end{remark}

Similar, stronger conditions on the cocycle was used by Conley and Tamuz \cite{ConleyTamuz}.
They designed a procedure for solving the \emph{unfriendly coloring} problem on graphs of degree bounded by $\Delta$, and showed that this procedure terminates off of a $\mu$-null set for every quasi-invariant Borel probability measure $\mu$ under the assumption that $1-1/\Delta\le\rho_\mu(x,y)\le 1+1/\Delta$ for every edge $\{x,y\}$ (this includes for example the case when $\mu$ is invariant).
Finding a measurable unfriendly coloring for a general Borel probability measure remains an interesting open problem.
In general, it would be nice to investigate if our condition has another applications in the context of local graph coloring problems.


The paper is structured as follows: in \cref{sec:preliminaries} we set the notation and recall basic results, in \cref{sec:EdgeCol} we define the augmenting chains that we consider in this paper, so-called $3$-step Vizing chains, and estimate how many of them can be attached to an uncolored edge, in \cref{sec:ImproveCol} we describe an infinite procedure that improves a given coloring so that it does not contain augmenting chains of small weight, in \cref{sec:Cocycle} we prove  \cref{thm:CocycleIntro}, that is, we show how to modify a given quasi-invariant measure to an equivalent measure that has bounded cocycle on edges, in \cref{sec:DoubleCount} we describe the double counting argument that estimates the number of uncolored edges, and finally in \cref{sec:Proof} we combine all the results to prove \cref{thm:MainIntro}.

\subsection*{Acknowledgements} 
	I would like to thank Oleg Pikhurko for many insightful discussions and constant support, to G\'{a}bor Elek for pointing out the reference \cite{KuhnCocycle}, and to the anonymous referee for valuable comments.
    The author was supported by Leverhulme Research Project Grant RPG-2018-424.

\section{Preliminaries}\label{sec:preliminaries}

Our basic descriptive set theory reference is \cite{kechrisclassical}, see also \cite{kechris_marks2016descriptive_comb_survey,pikhurko2021descriptive_comb_survey}.
We mostly follow the notation developed in \cite{grebik2020measurable}.
Let us point out that $\mathbb{N}$ contains $0$.
If $k\in \mathbb{N}\setminus \{0\}$, then we write $[k]=\{1,2,\dots,k\}$.
If $S$ is a set, then we write $[S]^k$ for the set of $k$-element subsets of $S$.
In particular, $[S]^2$ is the set of unordered pairs of $S$.
We follow a standard graph theoretic notation, i.e., a graph is a pair $G=(V,E)$, where $V\subseteq [E]^2$, $\deg_G(x)$ is the number of neighbors of $x$ in $G$, $\dist_G(x,y)$ denotes the graph distance of $x,y\in V$ in $G$, and $N_G(x)$ denotes the set of all edges of $G$ that are adjacent to $x$, i.e., $N_G(x)=\{e\in E:x\in e\}$.\\

{\bf Borel graphs \ \ }
A \emph{Borel graph} is a triplet $\fG=(V,\fB,E)$, where $(V,\fB)$ is a standard Borel space, $(V,E)$ is a graph and $E$ is a Borel subset of $[V]^2$ (endowed with the natural $\sigma$-algebra inherited from the product $\sigma$-algebra $\fB\times \fB$ on $V\times V$).
We say that $\fG$ is of \emph{bounded maximum degree} if there is $\Delta\in \mathbb{N}$ such that $\deg_\fG(x)\le \Delta$ for every $x\in V$.
We write $\Delta(\fG)$ for smallest such $\Delta$.
We denote the connectedness relation of $\fG$ as $F_\fG$.
That is, $F_\fG$ is an equivalence relation on $X$, and we write $[x]_\fG$ for the $\fG$-connectivity component of $x\in V$.
If $\Delta(\fG)<\infty$, then $F_\fG$ is a Borel subset of $V\times V$ and $|[x]_{\fG}|$ is at most countable.
Recall that an equivalence relation $F$ on a standard Borel space $(X,\fD)$ that is Borel as a subset of $X\times X$ and each $F$-equivalence class is at most countable is called a \emph{countable Borel equvialence relation (CBER)} on $X$, see \cite{KechrisCBER}.
In particular, if $\Delta(\fG)<\infty$, then $F_\fG$ is a CBER on $V$.
Given $A\subseteq X$, we write $[A]_F$ for the \emph{$F$-saturation of $A$ in $F$}, i.e., the smallest $F$-invariant subset of $X$ that contains $A$.
If $[A]_F=A$, then we say that $A$ is \emph{$F$-invariant}.
In case $F=F_\fG$ we write simply $[A]_\fG$, and say that $A$ is $\fG$-invariant if it is $F_\fG$-invariant.
Observe that if $A$ is a Borel set, then so is $[A]_F$.

The \emph{line graph} of a Borel graph $\fG=(V,\fB,E)$ is the Borel graph $\fE=(E,\fC,I_\fG)$, where $\fC$ is the $\sigma$-algebra on $E$ inherited form $[V]^2$ and $(E,I_\fG)$ is the line graph of $(V,E)$, i.e., $\{e,f\}\in I_\fG$ if and only if $e$ and $f$ share exactly one vertex.
Observe that if $\Delta(\fG)<\infty$, then $\Delta(\fE)\le 2\Delta(\fG)-2$.
Similarly as above, $F_\fE$ is the CBER on $E$ induced by the connectivity components of $\fE$.

A \emph{Borel chromatic number $\chi_\fB(\fG)$} of a Borel graph $\fG$ is the minimal $k\in \mathbb{N}$ such that there is a Borel proper vertex coloring $d:V\to [k]$ of $\fG$.
Note that the subscript $\chi_\fB$ refers to the corresponding $\sigma$-algebra.
A \emph{Borel chromatic index $\chi'_\fB(\fG)$} is defined as $\chi_\fC(\fE)$.
That is $\chi'_\fB(\fG)=k$ if there is a Borel proper vertex coloring $c:E\to [k]$ of $\fE$.\\

{\bf Measures \ \ }
The set of all \emph{Borel probability measures} on a standard Borel space $(X,\fB)$ is denoted as $\fP(X)$.
Let $F$ be a  CBER on $X$ and $\mu\in \fP(X)$.
We say that $\mu$ is \emph{$F$-quasi invariant} if $\mu([A]_F)=0$, whenever $\mu(A)=0$.
If $\fG=(V,\fB,E)$ is a Borel graph then we say that $\mu\in \fP(V)$ is \emph{$\fG$-quasi-invariant} if it is $F_\fG$-quasi-invariant

\begin{proposition}[Proposition~3.1 and 3.2 in \cite{grebik2020measurable}]\label{pr:BasicMeasure}
Let $\fG=(V,\fB,E)$ be a Borel graph such that $\Delta(\fG)<\infty$, $\fE=(E,\fC,I_\fG)$ be its line graph and $\mu\in \fP(V)$.
Then there is $\hat\mu\in \fP(E)$ that is $\fE$-quasi-invariant and satisfies
$$\mu\left(\left\{x\in V:\exists e\in A \ x\in e\right\}\right)=0$$
for every $A\subseteq E$ that satisfy $\hat\mu(A)=0$. 
\end{proposition}
\bpf
    By \cite[Proposition~3.2]{grebik2020measurable}, we find $\tilde \mu\in \fP(V)$ that is $\fG$-quasi invariant and satisfies $\mu([A]_\fG)=\tilde \mu([A]_\fG)$ for every $A\subseteq V$.
    By \cite[Proposition~3.1]{grebik2020measurable}, we find $\hat \mu\in \fP(E)$ that is $\fE$-quasi invariant and satisfies
    $$\tilde{\mu}\left(\left\{x\in V:\exists e\in A \ x\in e\right\}\right)<\Delta(\fG)\epsilon$$
    for every $A\subseteq E$ that satisfies $\hat \mu(A)<\epsilon$.
    In particular, if $\hat \mu(A)=0$, then
    \begin{equation*}
        \begin{split}
            \mu\left(\left\{x\in V:\exists e\in A \ x\in e\right\}\right)\le & \  \mu\left(\left[\left\{x\in V:\exists e\in A \ x\in e\right\}\right]_\fG\right)\\
            = & \ \tilde \mu\left(\left[\left\{x\in V:\exists e\in A \ x\in e\right\}\right]_\fG\right)=0
        \end{split}
    \end{equation*}
    as $\tilde{\mu}\left(\left\{x\in V:\exists e\in A \ x\in e\right\}\right)=0$.
\epf

A fundamental tool in the study of quasi-invariant measures is the Radon-Nikodym cocycle.
Let $(X,\fB)$ be a standard Borel space, $F$ be a CBER on $X$ and $\mu\in \fP(X)$ be $F$-quasi-invariant.
Then the \emph{Radon--Nikodym cocycle (of $\mu$ with respect to $F$)} is a Borel function $\rho_{\mu,F}:F\to \mathbb{R}_{>0}$ with the property that 
$$\mu(g(C))=\int_C \rho_{\mu,F}(x,g(x)) \ d\mu(x)$$
for every $C\in \fB$ and injective Borel map $g:C\to X$ such that $(x,g(x))\in F$.
It is a standard fact that the Radon--Nikodym cocycle exists and it is unique up to null-sets, that is, if $\rho$ and $\rho'$ are two Radon--Nikodym cocycles of $\mu$ with respect to $F$, then there is a $\mu$-conull $F$-invariant set $A\subseteq X$ such that
$$\rho\upharpoonright (F\cap (A\times A))=\rho'\upharpoonright (F\cap (A\times A)),$$
see \cite{topics}.
The following statement summarizes the properties of cocycles that we need.

\begin{proposition}[Chapter~II, Section~8 in \cite{topics}]\label{pr:BasicCocycle}
    Let $(X,\fB)$ be a standard Borel space, $F$ be a CBER on $X$ and $\mu\in \fP(X)$ be $F$-quasi-invariant.
    Then the Radon--Nikodym cocycle $\rho_{\mu,F}:F\to \mathbb{R}_{>0}$ satisfies:
    \begin{enumerate}
        \item there is a $\mu$-conull $F$-invariant set $A\in \fB$ such that $\rho_{\mu,F}(x,y)\rho_{\mu,F}(y,z)=\rho_{\mu,F}(x,z)$ for any $x,y,z\in A$ such that $y,z\in [x]_F$,
        \item (mass transport) we have
        $$\int_{X} \sum_{x\in [y]_F}{\bf F}(x,y) \ d \mu(y)=\int_{X} \sum_{y\in [x]_F} {\bf F}(x,y)\rho_{\mu,F}(x,y) \ d \mu(x)$$
        for any function ${\bf F}:F\to [0,\infty]$.
    \end{enumerate}
\end{proposition}

In case that $F=F_\fG$, we write $\rho_{\mu,\fG}$ instead of $\rho_{\mu,F_\fG}$.
If $F$ is clear from context, then we write simply $\rho_\mu$.

\begin{definition}
    Let  $\fG=(V,\fB,E)$ be a Borel graph such that $\Delta(\fG)<\infty$ and $\mu\in \fP(V)$ be $\fG$-quasi-invariant.
    We say that $\mu$ is \emph{$\fG$-bounded} if
    $$\frac{1}{4\Delta}\le \rho_\mu(x,y)\le 4\Delta$$
    for every $\{x,y\}\in E$.    
\end{definition}

We show in \cref{thm:MainCocycle} that every $\fG$-quasi-invariant $\mu\in \fP(V)$ is equivalent with a $\fG$-bounded measure $\nu\in \fP(V)$.

\section{Edge colorings}\label{sec:EdgeCol}

Let $\fG=(V,\fB,E)$ be a Borel graph such that $\Delta(\fG)<\infty$.
A \emph{partial (Borel proper edge) coloring} of $\fG$ is a partial Borel (with respect to the $\sigma$-algebra $\fC$ on $E$, in particular, $\dom(c)\in \fC$) map $c;E\to [\Delta(\fG)+1]$ that assigns different colors to different edges that share a vertex.
Usually we use lower case Greek letters for colors, e.g. $\alpha,\beta\in [\Delta(\fG)+1]$.
Given a partial coloring $c$, we define $m_c(x)$ to be the set of missing colors at $x\in V$, i.e., $m_c(x)=[\Delta(\fG)+1]\setminus \{c(e):e\in N_\fG(x)\}$.
We also write $U_c$ for the set of uncolored edges, i.e., $U_c=E\setminus \dom(c)$.

In order to improve a given partial coloring $c$, we utilize an idea from the proof of Vizing's theorem, so called \emph{Vizing chains}.
In general, given an uncolored edge $e\in U_c$, we want to find an injective augmenting sequence of edges $W_c(e)=(e_i)_{i\le k}$ such that $e_0=e$, $e_i\cap e_{i+1}\not=\emptyset$ for $i\le k$ and $e_i\not \in U_c$ for every $1\le i\le k$ with the property that keeping the colors outside of $W_{c}(e)$ intact but shifting the colors from $e_{i+1}$ to $e_i$ produces a different partial (proper) coloring $c'$ such that $m_{c'}(z_0)\cap m_{c'}(z_1)\not=\emptyset$, where $\{z_0,z_1\}=e_k$ is the last edge of $W_{c}(e)$.
Extending $c'$ by assigning any color from $m_{c'}(z_0)\cap m_{c'}(z_1)$ to $e_k$ then improves $c$ as we decreased the number of uncolored edges.
Observe that the difference between $c$ and $c'$ is contained in $W_{c}(e)$.

Various types of sequences $W_{c}(e)$ have been used in the literature.
In order to prove classical Vizing's theorem, one chooses $W_{c}(e)$ to be a concatenation of a so-called \emph{fan} and an \emph{alternating path}, also known as a \emph{Vizing chain}.
To prove an analogue of Vizing's theorem for graphings, the author and Pikhurko \cite{grebik2020measurable} iterated this process two times.
Namely, first we fix a Vizing chain, then truncate it at any edge on the alternating path, and then grow a second Vizing chain from that place; such a sequence is called an \emph{iterated Vizing chain}.
In order to devise efficient (both deterministic and randomized) local algorithms that produce a proper edge colorings with $\Delta+1$ colors on finite graphs, Bernshteyn \cite{BernshteynVizing} iterated this process $\log(n)$ times, where $n$ is the number of vertices of the graph.
Bernshteyn called the produced injective sequence a \emph{multi-step Vizing chain}.
In this paper, being ideologically more closer to \cite{grebik2020measurable}, we iterate the process $3$ times.
We follow the terminology of Bernshteyn and call this augmenting chain a \emph{$3$-step Vizing chain}.
The estimate for the number of $3$-step Vizing chains that can be assigned to a given uncolored edge follows the computation from \cite{grebik2020measurable}.
It seems plausible that one can get similar estimate by adapting the results from finite graphs, but we decided to follow the more direct path of generalizing \cite{grebik2020measurable}.

\subsection{$3$-step Vizing chains}
We recall the notation from \cite[Section 2]{grebik2020measurable} and refer the reader there for basic results about this notation.
Fix a Borel graph $\fG=(V,\fB,E)$ such that $\Delta(\fG)<\infty$ and a partial edge coloring $c;E\to [\Delta(\fG)+1]$.
Set $\Delta=\Delta(\fG)$ and fix an ordering of the set of colours $[\Delta+1]$.

A {\it chain} is a sequence $P=(e_0,\dots )$ of edges of $\fG$ such that for every index $i\in \mathbb{N}$ with $e_i,e_{i+1}$ being in $P$ we have $e_i\cap e_{i+1}\not=\emptyset$, that is, every two consecutive edges in $P$ intersect. 
Let $l(P)=|P|$ denote the \emph{length} of the chain $P$, i.e., the number of edges in $P$. 
Note that a chain can be finite (possibly empty) or infinite; thus $l(P)\in \mathbb{N}\cup\{\infty\}$ and, if $P$ is finite, then $P=(e_0,\dots, e_{l(P)-1})$. The convention of labeling the first edge as $e_0$ allows us to write $P=(e_i)_{i<l(P)}$, regardless of whether $P$ is finite or not.
If $l(P)=\infty$, then  we define $l(P)-1=\infty$ in order to avoid case by case statements in several places.

We call $e_{i-1}$ the \emph{$i$-th} edge of $P$.
For an edge $f$ that occurs exactly once in $P$, let its \emph{index} $i(f)$ be $i\ge 1$ such that $f=e_{i-1}$, that is, the index of the $i$-th edge is $i$. Also, for $i\le l(P)$, let $P_i=(e_j)_{j<i}$ denote the \emph{$i$-th prefix} of~$P$ (which consists of the first $i$ edges from $P$). We have, for example, that $P_{l(P)}=P$. 
For chains $P$ and $Q$, we write $P\sqsubseteq Q$ if $P=Q_{l(P)}$, that is, $P$ is a prefix of~$Q$.
If $P$ is a finite chain with the last edge $e$ and $Q$ is a chain with the first edge $f$ and $e\cap f\not=\emptyset$, then we write $P^\frown Q$ for the chain that is the concatenation of $P$ and~$Q$.

Let us call a chain $P=(e_i)_{i<l(P)}$ a \emph{path} if $P$ is empty, or if every vertex $z\in V$ belongs to at most 2 edges from $P$ and there is a vertex that belongs only to $e_0$. (In other words, $P$ is a finite path with a fixed direction on edges or an infinite one-sided ray, where no self-intersections are allowed.) Also, a chain  $P$ is called a \emph{cycle} if $P$ is non-empty and every vertex belongs to 0 or 2 edges of $P$. (These are just finite cycles, having some edge and direction fixed.) 

\begin{definition}[Definition~2.2 in \cite{grebik2020measurable}]\label{cAdmiss}
We say that a chain $P=(e_i)_{i<l(P)}$ is 
\begin{enumerate}
\item \emph{edge injective} if every edge appears at most once in $P$, that is, for every $0\le i<j< l(P)$ we have that $e_i\not=e_j$,

\item \emph{$c$-shiftable} if $l(P)\ge 1$, $P$ is edge injective,  $e_0\in U_c$
and $e_j\in \dom(c)$ for every $1\le j<l(P)$ (that is, if $P$ is non-empty with no edge repeated and $e_0$ is the unique uncoloured edge of $P$);

\item \emph{$c$-proper-shiftable} if $P$ is $c$-shiftable and $c_P;E\to [\Delta+1]$ is a partial coloring, where $c_P$ is \emph{the shift of $c$ along $P$} (or \emph{$P$-shift of $c$} for short)
which is defined as 
   \begin{itemize}
   	\item $\dom(c_P)=\dom(c)\cup\{e_0\}\setminus \{e_{l(P)-1}\}$ where we put $\{e_{l(P)-1}\}=\emptyset$ if $l(P)=\infty$,
        \item $c_P(e_i)=c(e_{i+1})$ for every $i+1<l(P)$,
        \item $c_P(f)=c(f)$ for every $f\in \dom(c)\setminus P$;
    \end{itemize}
   \item \emph{$c$-augmenting} if $P$ is $c$-proper-shiftable and either $l(P)=\infty$ or $P$ is finite with $m_{c_P}(x)\cap m_{c_P}(y)\not=\emptyset$ where $x\not=y$ are the vertices of the last edge $e_{l(P)-1}$ of~$P$.
\end{enumerate}
\end{definition}

Next we describe the building blocks that will be used to build $3$-step Vizing chains:\\

{\bf Alternating path}.
Let $x\in V$ and $\alpha,\beta\in [\Delta+1]$ be different colours such that $\beta\in m_c(x)$.
Then there is a unique maximal chain $P=(e_i)_{i<l(P)}$ such that $x\in e_0$ if $l(P)>0$, $x\not\in e_1$ if $l(P)>1$, and $c(e_{i})=\alpha$ (resp.\ $c(e_{i})=\beta$) for every $i<l(P)$ that is even (resp.\ odd).
We call this unique maximal chain {\it the (alternating) $\alpha/\beta$-path starting at $x\in V$} and denote it as $P_c(x,\alpha/\beta)$. 
If $P_c(x,\alpha/\beta)$ is finite and non-empty, then we call the unique $y\in V$ such that $|\{f\in P_c(x,\alpha/\beta):y\in f\}|=1$ and $y\not=x$ the {\it last vertex of $P_c(x,\alpha/\beta)$}.
If $P_c(x,\alpha/\beta)$ is empty (which happens exactly when $\alpha\in m_c(x)$), then the \emph{last vertex} is~$x$.
Whenever we write $P_c(x,\alpha/\beta)$ we always assume that the condition that $\beta\in m_c(x)$ is satisfied. 
Observe that the colors on the chain alternate between $\alpha$ and $\beta$ (starting with $\alpha$) and we never return to a vertex we have previously visited (and thus the edges in $P$ form a path). \\

{\bf Fan}.
Let $e\in U_c$ and $x\in e$.
We define the {\it maximal fan around $x$ starting at $e$}, in symbols $F_c(x,e)$, as a (finite) chain $P=(e_0,e_1,\dots,e_k)$ such that $x\in e_i$ for every $i\le k$ and if we denote the other vertex in $e_i$ by $v_i$ then the following statements are satisfied
\begin{enumerate}
    \item $e_0=e$,
    \item $P$ is edge injective,
    \item $c(e_{i+1})\in m_c(v_i)$ for every $i<k$ and $c(e_{i+1})$ is the minimal color available in the $i$-th step,
    where we say that a color $\alpha$ is \emph{available in the $i$-th step} if $\alpha\in m_c(v_i)$,
    \item $(e_0,\dots,e_k)$ is maximal with these properties.
\end{enumerate}
 
{\bf Conditional fan}.
We generalize, but only formally, the following definition from \cite[Section~2.5]{grebik2020measurable}.
Take $P_c(x,\alpha/\beta)$ for some $x\in V$.
Let $f\in P_c(x,\alpha/\beta)$ be a an edge that is not first nor last in $P_c(x,\alpha/\beta)$ and $y\in V$ be the last vertex of $P_c(x,e)_{i(f)}$.
We define {\it the maximal $\alpha/\beta$-conditional fan starting at $f$}, denoted as $F_c(x,\alpha/\beta,y)$, as a chain $P=(g_0,\dots, g_m)$ such that $y\in g_i$ for every $i\le m$ and, if we denote the other vertex of $g_i$ by $u_i$, then the following is satisfied
\begin{enumerate}
    \item\label{it:cond1} $g_0=f$,
    \item\label{it:cond2} $P$ is edge injective,
    \item\label{it:cond3} $c(g_{i+1})\in m_c(u_{i})$ and it is the minimal available color,
    \item \label{it:cond4} $\alpha,\beta \not \in m_c(u_i)$ for every $i<m$,
    \item\label{it:cond5}  if $\alpha,\beta\not\in m_c(u_m)$, then $(g_0,\dots, g_m)$ is maximal with the properties above.
\end{enumerate}
Note that we should rather write $u_i^f$, $g_i^f$ and $y^f$ to stress that those objects depend on the choice of $f$.
This will be however omitted in the cases when we work with only one $f$.\\

Now we are ready to for the main definition.

\begin{definition}[$3$-step Vizing chain]\label{def:VizingChain}
    Let $\fG=(V,\fB,E)$ be a Borel graph such that $\Delta(\fG)<\infty$, $c;E\to [\Delta(\fG)+1]$ be a partial edge coloring and $e\in U_c$.
    We say that a $c$-augmenting chain $W_c(e)=(e_i)_{i< l(W_c(e))}$, where $l(W_c(e))\in \mathbb{N}\cup \{\infty\}$, is a \emph{$3$-step Vizing chain (at $e$)} if there are pairwise different vertices $y_1,y_2,y_3,z_1,z_2,z_3\in V$, and colors $\alpha_i,\beta_i\in [\Delta(\fG)+1]$ for $i\in \{1,2,3\}$ such that
    $$W_c(e)=\left(F^1_c\right)^\frown \left(P^1_c\right)^\frown\left(F^2_c\right)^\frown\left(P^2_c\right)^\frown\left(F^3_c\right)^\frown\left(P^3_c\right),$$
    where
            \begin{enumerate}
                \item $F^1_c\sqsubseteq F_c(y_1,e)$, $F^2_c\sqsubseteq F_c(z_1,\alpha_1/\beta_1,y_2)$ and $F^3_c\sqsubseteq F_c(z_2,\alpha_2/\beta_2,y_3)$,
                \item $P^i_c\sqsubseteq P_c(z_i,\alpha_i/\beta_i)$ for every $i\in \{1,2,3\}$,
                \item if $T\in \{F^i_c\}_{i=1}^3\cup \{P^i_c\}_{i=1}^3$ satisfies $T=\emptyset$, then every $S$ to the right from $T$ in the definition of $W_c(e)$ is empty as well, i.e., $W_c(e)$ is built by \emph{at most} three iterations of the ``Vizing Chain'' construction,
            \end{enumerate}
\end{definition}

\subsection{Construction of $3$-step Vizing chains}
Let $e\in U_c$ and $x\in e$ be fixed.
First, we describe a process that produces many chains of the form 
$$W_c(e)=\left(F^1_c\right)^\frown \left(P^1_c\right)^\frown\left(F^2_c\right)^\frown\left(P^2_c\right)^\frown\left(F^3_c\right)^\frown\left(P^3_c\right)$$
that satisfies (1)--(3) in \cref{def:VizingChain}, then we investigate how many of these chains are $3$-step Vizing chains, i.e., which ones are $c$-augmenting.
We handle each iteration separately.\\

{\bf Iteration I}.
We start with \cite[Section~2.4]{grebik2020measurable}.
Recall that the \emph{Vizing chain} $V_c(x,e)$ either consists of the fan $F_c(x,e)$, in case it is augmenting, or we have
$$V_c(x,e)={F_c(x,e)_{i+1}}^\frown P_c(v_i, \alpha/\beta),$$
where $i<l(F_c(x,e))$ (so-called first critical index) and $\{x,v_{i}\}$ is the last edge of the truncated fan $F_c(x,e)_{i_1+1}$.

\begin{claim}[Proposition~2.9 in \cite{grebik2020measurable}]\label{cl:GP1}
    If $F_c(x,e)$ is not $c$-augmenting, then there is $i<l(F_c(x,e))$ such that ${F_c(x,e)_{i+1}}^\frown P_c(v_i, \alpha/\beta)$ is $c$-augmenting, where $(x,v_i)$ is the last edge of the truncated fan $F_c(x,e)_{i+1}$ and $\beta$ is the smallest missing color at $v_i$.
    In particular, the Vizing chain $V_c(x,e)$ is a $3$-step Vizing chain.
    
    Moreover, if nonempty, then $P_c(v_{i}, \alpha/\beta)$ does not use the vertex $x$.
\end{claim}

Suppose that $P_c(v_i, \alpha/\beta)$ is non-empty.
Define $F^1_c=F_c(x,e)_{{i_1}+1}$, $y_1=x$, $z_1=v_{i}$ and $\alpha_1=\alpha$, $\beta_1=\beta$.
This concludes the first iteration of the construction.\\

{\bf Iteration II}.
Recall that $f^1\in P_c(z_1, \alpha_1/\beta_1)$ is \emph{suitable} \cite[Definition~2.10]{grebik2020measurable} if it is of graph distance $\ge 3$ from $F^1_c$, it is not the last edge of $P_c(z_1, \alpha_1/\beta_1)$ and $c(f^1)=\alpha_1$ (we remark that the last condition only helps with the notation and is otherwise irrelevant).
Let $y_2$ be the last vertex of $P_c(z_1, \alpha_1/\beta_1)_{i(f^1)}$ and set $F_c(x,e\leadsto f^1)=F_c(z_1,\alpha_1/\beta_1,y_2)$, where $F_c(z_1,\alpha_1/\beta_1,y_2)$ is the maximal $\alpha_1/\beta_1$-conditional fan starting at $f^1$.

\begin{claim}[Proposition~2.12 in \cite{grebik2020measurable}]
    Let $f^1 \in  P_c(z_1, \alpha/\beta)$ be suitable.
    Then $$P = {V_c(x,e)_{i(f^1) - 1}}^\frown F_c(x,e\leadsto f^1)$$ is $c$-proper-shiftable.
\end{claim}

According to the \emph{type} of $f^1$, as defined in \cite[Section~2.5]{grebik2020measurable}, it is possible to assign an injective sequence of edges $Q$ such that either $Q=F_c(x,e\leadsto f^1)$, or
$$Q={F_c(x,e\leadsto f^1)_{m+1}}^\frown P_c(u_m, \gamma/\delta),$$
where $m<l(F_c(x,e\leadsto f^1))$ (the so called second critical index), $\{y_2,u_m\}$ is the last edge of the truncated fan $F_c(x,e\leadsto f^1)_{m+1}$ and $\{\gamma,\delta\}$ is either equal to or disjoint from $\{\alpha_1,\beta_1\}$.
Using the technical notion of \emph{superb} edges, the following, again adapted to our terminology, is shown in \cite{grebik2020measurable}.

\begin{claim}[Proposition~2.15 and 2.17 in \cite{grebik2020measurable}]\label{cl:GP2}
    Let $f^1\in P_c(z_1, \alpha_1/\beta_1)$ be suitable and superb.
    The chain
    $$(F^1_c)^\frown (P_c(z_1, \alpha_1/\beta_1)_{i(f^1)-1})^\frown Q$$
    is a $3$-step Vizing chain.
\end{claim}

Suppose that $P_c(u_m, \gamma/\delta)$ is non-empty.
Define $P^1_c=P_c(z_1, \alpha_1/\beta_1)_{i(f^1)-1}$, $F^2_c=F_c(x,e\leadsto f^1)_{m+1}$, $z_2=u_m$ and $\alpha_2=\gamma$, $\beta_2=\delta$.
This concludes the second iteration of the construction.\\

{\bf Iteration III}.
We say that $f^2\in P_c(z_2,\alpha_2/\beta_2)$ is \emph{$2$-suitable} if
\begin{itemize}
    \item the graph distance of $f$ and $f^2$ is at least $3$ for every $f\in (F^1_c)^\frown (P^1_c)^\frown (F^2_c)$,
    \item it is not the last edge of $P_c(z_2,\alpha_2/\beta_2)$
    \item $c(f^2)=\alpha_2$.
\end{itemize}
Let $y_3$ be the last vertex of $P_c(z_2, \alpha_2/\beta_2)_{i(f^2)}$ and set
$$F_c(x,e\leadsto f^1\leadsto f^2)=F_c(z_2,\alpha_2/\beta_2,y_3),$$
where $F_c(z_2,\alpha_2/\beta_2,y_3)$ is the maximal $\alpha_2/\beta_2$-conditional fan starting at $f^2$.

\begin{proposition}\label{pr:TwoSuitable}
    Let $f^2 \in  P_c(z_2,\alpha_2/\beta_2)$ be $2$-suitable.
    Then
    $$P = (F^1_c)^\frown (P^1_c)^\frown (F^2_c)^\frown (P_c(z_2,\alpha_2/\beta_2)_{i(f^2)-1})^\frown F_c(x,e\leadsto f^1\leadsto f^2)$$
    is $c$-proper-shiftable.
\end{proposition}
\bpf
    Suppose that $c_P$ is not a partial coloring.
    By the definition, we find a vertex $v\in V$ such that $c_P\upharpoonright N_\fG(v)$ is not proper.
    By the fact that $f^2$ is $2$-suitable together with \cite[Proposition~2.12]{grebik2020measurable}, we see that $v\not\in f$ for any $f\in (F^1_c)^\frown (P^1_c)^\frown (F^2_c)$.
    Moreover, we must have $v\in f$ for some $f\in F_c(x,e\leadsto f^1\leadsto f^2)$ as the colors used around vertices that lie only on the path $P_c(z_2,\alpha_2/\beta_2)_{i(f^2)-1}$ do not change.
    Now the definition of $F_c(x,e\leadsto f^1\leadsto f^2)$, as the maximal $\alpha_2/\beta_2$-conditional fan starting at $f^2$, shows that no such $v\in V$ can exist, the argument is literally the same as in \cite[Proposition~2.12]{grebik2020measurable}. 
\epf

Suppose that $f^2$ is $2$-suitable and set $P^2_c=P_c(z_2,\alpha_2/\beta_2)_{i(f^2)-1}$.
Next we define various types of $2$-suitable edges and the notion of an amazing edge.
This is inspired by similar notions in \cite[Section~2.5]{grebik2020measurable}.

\vspace{+0.3cm}

We say that a $2$-suitable edge $f^2\in P_c(z_2,\alpha_2/\beta_2)$ is of \emph{type (a)} if 
$$P = (F^1_c)^\frown (P^1_c)^\frown (F^2_c)^\frown (P^2_c)^\frown F_c(x,e\leadsto f^1\leadsto f^2)$$
is $c$-augmenting.
Every edge of type (a) is said to be \emph{amazing}.
Define $F^3_c=F_c(x,e\leadsto f^1\leadsto f^2)$ and set
$$W_c(e,f^1,f^2)=(F^1_c)^\frown (P^1_c)^\frown (F^2_c)^\frown (P^2_c)^\frown (F^3_c).$$
Then the following is immediate from the definitions.
\begin{proposition}
    Let $f^2$ be of type (a).
    Then $W_c(e,f^1,f^2)$ is a $3$-step Vizing chain.
\end{proposition}

\vspace{+0.3cm}

We say that a $2$-suitable edge $f^2\in P_c(z_2,\alpha_2/\beta_2)$ is of \emph{type (b)} if it is not of type (a) and in the construction of the conditional fan $F_c(x,e\leadsto f^1\leadsto f^2)$ we encountered $\alpha_2$ or $\beta_2$.
Observe that as $f^2$ is not of type (a), we must have $\beta_2\in m_c(w_n)$, where $\{y_3,w_n\}$ is the last edge in $F_c(x,e\leadsto f^1\leadsto f^2)$.
Following the previous notation we say that $n$ is the \emph{third critical index}.
Define $F^3_c=F_c(x,e\leadsto f^1\leadsto f^2)$, $y_3=w_n$, $\alpha_3=\alpha_2$, $\beta_3=\beta_2$ and $P^3_c=P_c(z_3,\alpha_3/\beta_3)$.
We say that $f^2$ of type (b) is \emph{amazing}, if 
\begin{itemize}
    \item $P_{c_{Q}}(z_3,\alpha_3/\beta_3)=P_c(z_3,\alpha_3/\beta_3)$, where $Q=(F^1_c)^\frown (P^1_c)^\frown (F^2_c)^\frown (P^2_c)^\frown (F^3_c)$.
\end{itemize}

\begin{proposition}
    Let $f^2$ be of type (b) and amazing.
    Then
    $$W_c(e,f^1,f^2)=(F^1_c)^\frown (P^1_c)^\frown (F^2_c)^\frown (P^2_c)^\frown (F^3_c)^\frown (P^3_c)$$
    is a $3$-step Vizing chain.
\end{proposition}
\bpf
    It follows directly from the definitions that $W_c(e,f^1,f^2)$ satisfies (1)--(3) in \cref{def:VizingChain}.
    It remains to show that $W_c(e,f^1,f^2)$ is $c$-augmenting.
    This can be done by the same argument as in \cite[Proposition~2.15]{grebik2020measurable}.
    Namely, first observe that $c_Q$, where 
    $$Q=(F^1_c)^\frown (P^1_c)^\frown (F^2_c)^\frown (P^2_c)^\frown (F^3_c),$$
    is a proper coloring by the same reasoning as in \cref{pr:TwoSuitable}.
    Moreover, by the definition of $2$-suitable edge, we must have $\alpha_2\in m_{c_{Q}}(y_3)$.
    As $P_{c_{Q}}(z_3,\alpha_3/\beta_3)=P_c(z_3,\alpha_3/\beta_3)$, we have that $W_c(e,f^1,f^2)$ is edge injective and $\beta_2\in m_{c_{Q}(z_3)}$.
    This shows that $\{y_3,z_3\}^\frown P^3_c$ is $c_Q$-augmenting as $y_3$ cannot be the last vertex of $P^3_c$ (if it were, then $P_{c_{Q}}(z_3,\alpha_3/\beta_3)\not =P_c(z_3,\alpha_3/\beta_3)$ as $\alpha_2,\beta_2\not \in m_c(y_3)$).
    Hence $W_c(e,f^1,f^2)$ is $c$-augmenting as desired.
\epf

\vspace{+0.3cm}

We say that a $2$-suitable edge $f^2\in P_c(z_2,\alpha_2/\beta_2)$ is of \emph{type (c)} if it is not of type (a), or (b).
Let $\gamma$ be the smallest colour in $m_c(y_3)$.
The reason why we cannot extend $F_c(x,e\leadsto f^1\leadsto f^2)$ is the same as when we build the standard Vizing chain, see \cite[Proposition~2.8]{grebik2020measurable}, or \cite[Section2.5: Type II edge]{grebik2020measurable}.
Namely there is a colour $\delta$ and index
$$j<n=l(F_c(x,e\leadsto f^1\leadsto f^2))-1$$
such that $\delta$ is the minimal colour available in both $m_c(w_j)$ and $m_c(w_n)$.
It is clear that $\gamma\not=\delta$ because $f^2$ is not of Type (a) and $\{\alpha_2,\beta_2\}\cap \{\gamma,\delta\}=\emptyset$ because $f$ is not of Type (b).

Consider now the alternating $\gamma/\delta$-paths
$P_c(w_j,\gamma/\delta)$ and
$P_c(w_n,\gamma/\delta)$.
Our aim is to choose one of them, call it $Q$, and then define
$$W_c(e,f^1,f^2)=(F^1_c)^\frown (P^1_c)^\frown (F^2_c)^\frown (P^2_c)^\frown(F_c(x,e\leadsto f^1\leadsto f^2)_{\ell+1})^\frown Q,$$
where $\ell\in \{i,m\}$, depending on the choice of $Q$, is such that $W_c(x,f^1,f^2)$ is $c$-augmenting.
As in the case of type (b), we need to rule out some edges.
We say that a $2$-suitable $f^2\in P_c(z_2,\alpha_2/\beta_2)$ of type (c) is \emph{amazing} if, in the above notation, both of the following equalities hold 
\begin{itemize}
    \item $P_c(w_j,\gamma/\delta)=P_{c_R}(w_j,\gamma/\delta)$,
    \item $P_c(w_n,\gamma/\delta)=P_{c_R}(w_j,\gamma/\delta)$,
\end{itemize}
where $R=(F^1_c)^\frown (P^1_c)^\frown (F^2_c)^\frown (P^2_c)$.

Let $f^2\in P_c(z_2,\alpha_2/\beta_2)$ be of type (c) and \emph{amazing}.
We take $\ell\in \{j,n\}$ to be the index for which there is no $h\in P_c(u_j,\delta/\epsilon)$ such that $y_3\in h$, see \cite[Proposition~2.9]{grebik2020measurable}.
If both indices $j$ and $n$ satisfy this, then we put $\ell=j$ for definiteness.
We call this index $\ell$ the {\it third critical index}.
We define $F^3_c=F_c(x,e\leadsto f^1\leadsto f^2)_{\ell+1}$, $z_3=w_\ell$, $\alpha_3=\gamma$, $\beta_3=\delta$ and $P^3_c=P_c(w_\ell,\alpha_3/\beta_3)$.

\begin{proposition}
    Let $f^2$ be of type (c) and amazing.
    Then
    $$W_c(e,f^1,f^2)=(P^1_c)^\frown (F^2_c)^\frown (P^2_c)^\frown (F^3_c)^\frown (P^3_c)$$
    is a $3$-step Vizing chain.
\end{proposition}
\bpf
    It follows directly from the definitions that $W_c(e,f^1,f^2)$ satisfies (1)--(3) in \cref{def:VizingChain}.
    It remains to show that $W_c(e,f^1,f^2)$ is $c$-augmenting.
    This can be done by the same argument as in \cite[Proposition~2.17]{grebik2020measurable}.
    Namely, first observe that $c_Q$, where 
    $$Q=(F^1_c)^\frown (P^1_c)^\frown (F^2_c)^\frown (P^2_c)^\frown (F^3_c),$$
    is $c$-proper shiftable by the same reasoning as in \cref{pr:TwoSuitable}.
    In particular, $Q$ is edge injective.
    As $y_3\not \in f$ for any $f\in P^3_c$ by the definition of the third critical index and $P_{c_{R}}(z_3,\alpha_3/\beta_3)=P_c(z_3,\alpha_3/\beta_3)$, we infer that $W_c(e,f^1,f^2)$ is edge injective.
    Similar argument shows that $P_{c_{Q}}(z_3,\alpha_3/\beta_3)=P_{c_{R}}(z_3,\alpha_3/\beta_3)=P_c(z_3,\alpha_3/\beta_3)$.
    Moreover, by the definition of $2$-suitable edge, we must have $\alpha_3\in m_{c_{Q}}(y_3)$.
    This shows that $\{y_3,z_3\}^\frown P^3_c$ is $c_Q$-augmenting.
    Hence $W_c(e,f^1,f^2)$ is $c$-augmenting as desired.
\epf

Altogether, we just proved the following statement.

\begin{theorem}\label{thm:MainChains}
    Let $e\in U_c$, $x\in e$, $f^1$ be superb and $f^2$ be amazing as defined above.
    Then $W_c(e,f^1,f^2)$ is a $3$-step Vizing chain.
\end{theorem}

\subsection{How many $3$-step Vizing chains are there}

Let $e\in U_c$ and $x\in e$.
We say that $e$ is \emph{$K$-bad for $c$}, where $K\in \mathbb{N}$, if every $3$-step Vizing chain $W_c(e)$ at $e$ satisfies $l(W_c(e)))\ge 2K+2\Delta$.
In the following claims we use the notation from previous section.

\begin{proposition}
    Let $K\in \mathbb{N}$, $e\in U_c$ be $K$-bad for $c$ and $x\in e$.
    Then
    \begin{itemize}
        \item $l(P_c(v_i,\alpha/\beta))\ge K$, where $P_c(v_i,\alpha/\beta)$ is the alternating path from the first iteration,
        \item $l(P_c(u_m,\gamma/\delta))\ge \frac{K}{2}$ for every superb edge $f^1$ such that $i(f^1)\le \frac{K}{2}$, where $P_c(u_m,\gamma/\delta)$ is the alternating path in the second iteration that corresponds to $f^1$.
    \end{itemize}
\end{proposition}
\bpf
    Both chains from \cref{cl:GP1} and \cref{cl:GP2} are $3$-step Vizing chains.
    As $e$ is $K$-bad for $c$ and both chains contain at most two fans that each contain at most $\Delta$ edges, we conclude that $l(P_c(v_i,\alpha/\beta))\ge K$ and $l(P_c(u_m,\gamma/\delta))\ge \frac{K}{2}$ under the assumption that $f^1$ is superb and satisfies $i(f^1)\le \frac{K}{2}$.
\epf

\begin{claim}[Proposition~2.20 in \cite{grebik2020measurable}]\label{cl:GP3}
    Let $K\in \mathbb{N}$, $e\in U_c$ be $K$-bad for $c$ and $x\in e$.
    Then there are colors $\{\gamma,\delta\}\subseteq [\Delta+1]$ and at least
    $$\frac{1}{3(\Delta+1)^2}\left(\frac{K}{4}-\Delta^5-1\right)-2\Delta^3$$
    many superb edges $f^1\in P_c(v_i,\alpha/\beta)$ such that $i(f^1)\le \frac{K}{2}$ and the alternating path in the second iteration that corresponds to $f^1$ is a $\gamma/\delta$-path.
\end{claim}

\begin{definition}\label{def:NumberPairs}
    Let $K\in \mathbb{N}$, $e\in U_c$ be $K$-bad for $c$, $x\in e$ and $\{\gamma,\delta\}\subseteq [\Delta+1]$ be as in \cref{cl:GP3}.
    Define $\fV_c(e)$ to be the set of pairs $(f^1,f^2)$ such that $f^1\in P_c(v_i,\alpha/\beta)$ is superb and satisfies $i(f^1)\le \frac{K}{2}$ and $f^2\in P_c(u_m,\gamma/\delta)$ is amazing and satisfies $i(f^2)\le \frac{K}{2}$ (where the index $i(f^2)$ is taken with respect to $P_c(u_m,\gamma/\delta)$), where $P_c(u_m,\gamma/\delta)$ is the alternating path in the second iteration that corresponds to $f^1$.
\end{definition}

Our aim is to estimate the cardinality of $\fV_c(e)$ as this gives a lower bound on the cardinality of the set of all $3$-step Vizing chains at $e$.
In fact, we bound the cardinality of the projection of $\fV_c(e)$ to the second coordinate as this clearly gives a lower bound for the cardinality of $\fV_c(e)$.

The following proposition gives a sufficient condition for an edge to be amazing.

\begin{proposition}\label{pr:AmazingEdge}
    Let $K\in \mathbb{N}$, $e\in U_c$ be $K$-bad for $c$, $x\in e$ and $\{\gamma,\delta\}\subseteq [\Delta+1]$ be as in \cref{cl:GP3}.
    Suppose that $f^1\in P_c(v_i,\alpha/\beta)$ is superb, $i(f^1)\le \frac{K}{2}$ and pick $f^2$ on the the alternating path $P_c(u_m,\gamma/\delta)$ (that corresponds to $f^1$ in the second iteration).
    Assume that
    \begin{enumerate}
        \item $f^2 $ is $2$-suitable, 
        \item there is no alternating path $P_c(w,\iota/\kappa)$ such that, simultaneously, $w\in V$ is of $\fG$-distance $1$ from $y_3\in f^2$ and $P_c(w,\iota/\kappa)$ is of distance at most $3$ from $(F^1_c)^\frown (P^1_c)^\frown(F^2_c)$, where $f^1$ is the first edge of $F^2_c$.
    \end{enumerate}
    Then $f^2$ is amazing.
    In particular, $(f^1,f^2)\in \fV_c(e)$.
\end{proposition}
\bpf
    Suppose that the conditions are satisfied.
    If $f^2$ is of type (a), then it is amazing.
    
    If $f^2$ is of type (b), then we need to verify that $P_{c_{Q}}(z_3,\alpha_3/\beta_3)=P_c(z_3,\alpha_3/\beta_3)$, where $Q=(F^1_c)^\frown (P^1_c)^\frown (F^2_c)^\frown (P^2_c)^\frown (F^3_c)$.
    As $\dist_\fG(y_3,z_3)=1$, we have that $P_c(z_3,\alpha_3/\beta_3)$ avoids $(F^1_c)^\frown (P^1_c)^\frown (F^2_c)$.
    Hence if $P_{c_{Q}}(z_3,\alpha_3/\beta_3)\not =P_c(z_3,\alpha_3/\beta_3)$ it must be the case that $y_3$ is covered by $P_c(z_3,\alpha_3/\beta_3)$.
    This can only happen if $P_c(z_3,\alpha_3/\beta_3)$ contains $P^2_c$ as $\alpha_3=\alpha_2$ and $\beta_3=\beta_2$.
    But that is not possible as $P^2_c$ is of distance $1$ from $(F^1_c)^\frown (P^1_c)^\frown (F^2_c)$.

    If $f^2$ is of type (c), then we need to verify that $P_c(w_j,\gamma/\delta)=P_{c_R}(w_j,\gamma/\delta)$ and $P_c(w_n,\gamma/\delta)=P_{c_R}(w_j,\gamma/\delta)$.
    This follows easily as $\{\gamma,\delta\}\cap \{\alpha_2,\beta_2\}=\emptyset$ and both $P_c(w_j,\gamma/\delta)$, $P_c(w_n,\gamma/\delta)$ avoid $(F^1_c)^\frown (P^1_c)^\frown (F^2_c)$.
\epf

\begin{proposition}\label{pr:LowerBound}
    Let $K\in \mathbb{N}$, $e\in U_c$ be $K$-bad for $c$, $x\in e$ and $\{\gamma,\delta\}\subseteq [\Delta+1]$ be as in \cref{cl:GP3}.
    Define $f^2\in \mathcal{S}_c(e)$ if $c(f^2)=\gamma$ and there is $f^1\in P_c(v_i,\alpha/\beta)$ such that $f^1$ is superb, $i(f^1)\le \frac{K}{2}$, $f^2 \in P_c(u_m,\gamma/\delta)$ and $i(f^2)\le\frac{K}{2}$ (the index is taken with respect to $P_c(u_m,\gamma/\delta)$ and $P_c(u_m,\gamma/\delta)$ corresponds to $f^1$ in the second iteration), and at least one of the items from \cref{pr:AmazingEdge} is not satisfied.
    Then we have
    $$|\mathcal{S}_c(e)|\le 4(\Delta+1)^4\sum_{r=0}^4 (2\Delta)^r\left(\frac{K}{2}+\Delta\right)$$
\end{proposition}
\bpf
    Let $f^1$ and $f^2$ are as above and assume that item (1) from \cref{pr:AmazingEdge} is not satisfied.
    As $c(f^2)=\gamma$, we know that either $f^2$ is the last edge on $P_c(u_m,\gamma/\delta)$ or it is of distance at most $4$ from $(F^1_c)^\frown P_c(v_i,\alpha/\beta)_{\frac{K}{2}+1}$.
    There are at most $\frac{K}{2}$ many edges that satisfy the former condition, and we have 
    \begin{equation}\label{eq:BadEdges}
        \left|\left\{f\in E:\dist_\fE(f,(F^1_c)^\frown P_c(v_i,\alpha/\beta)_{\frac{K}{2}+1})\le 4\right\}\right|\le\sum_{r=0}^4 (2\Delta)^r\left(\frac{K}{2}+\Delta\right)
    \end{equation}
    which gives upper bound on the latter condition.

    Suppose that (2) from \cref{pr:AmazingEdge} is not satisfied.
    That means that there is a path $P_c(w,\iota/\kappa)$, where $w$ is of distance one from $y_3$, that intersect
    $$B=\left\{f\in E:\dist_\fE(f,(F^1_c)^\frown P_c(v_i,\alpha/\beta)_{\frac{K}{2}+1})\le 4\right\}.$$
    Every edge from $B$ is an element of $2(\Delta+1)^2$ many such paths.
    Together with \eqref{eq:BadEdges}, we conclude that there are at most
    $$2(\Delta+1)^4\sum_{r=0}^4 (2\Delta)^r\left(\frac{K}{2}+\Delta\right)$$
    edges $f^2$ that do not satisfy (2) from \cref{pr:AmazingEdge}.

    Summing these three bounds gives the desired estimate.
\epf

\begin{proposition}\label{pr:UpperBound}
    Let $K\in \mathbb{N}$, $e\in U_c$ be $K$-bad for $c$, $x\in e$ and $\{\gamma,\delta\}\subseteq [\Delta+1]$ be as in \cref{cl:GP3}.
    Define $f^2\in \mathcal{T}_c(e)$ if $c(f^2)=\gamma$ and there is $f^1\in P_c(v_i,\alpha/\beta)$ such that $f^1$ is superb, $i(f^1)\le \frac{K}{2}$, $f^2 \in P_c(u_m,\gamma/\delta)$ and $i(f^2)\le\frac{K}{2}$ (the index is taken with respect to $P_c(u_m,\gamma/\delta)$ and $P_c(u_m,\gamma/\delta)$ corresponds to $f^1$ in the second iteration).
    Then we have
    $$|\mathcal{T}_c(e)|\ge \frac{K}{4\Delta^2} \left(\frac{1}{3(\Delta+1)^2}\left(\frac{K}{4}-\Delta^5-1\right)-2\Delta^3\right)$$
\end{proposition}
\bpf
    By \cref{cl:GP3}, there is $\{\gamma,\delta\}\subseteq [\Delta+1]$ and at least 
    $$\frac{1}{3(\Delta+1)^2}\left(\frac{K}{4}-\Delta^5-1\right)-2\Delta^3$$
    many superb edges $f^1\in P_c(v_i,\alpha/\beta)$ such that $i(f^1)\le \frac{K}{2}$ and the alternating path in the second iteration that corresponds to $f^1$ is a $\gamma/\delta$-path.
    For each such $f^1$, there are at least $\frac{K}{4}$ edges of the corresponding $P_c(u_m,\gamma/\delta)$ that have color $\gamma$.
    This shows that there are at least 
    $$\frac{K}{4}\left(\frac{1}{3(\Delta+1)^2}\left(\frac{K}{4}-\Delta^5-1\right)-2\Delta^3\right)$$
    pairs $(f^1,f^2)$ that satisfy the conditions above.
    To estimate the size of $\mathcal{T}_c(e)$ we need to compute the number of pairs $(f^1,f^2)$ to which a given edge $f^2$ contributes.

    Every $f^2$ can reach $f_1$ by following the $\gamma/\delta$ path in one of its two directions, and then there are $\Delta^2$ many choices for $f^1$.
    Altogether,
    $$|\mathcal{T}_c(e)|\ge\frac{K}{4\Delta^2} \left(\frac{1}{3(\Delta+1)^2}\left(\frac{K}{4}-\Delta^5-1\right)-2\Delta^3\right)$$
    as needed.
\epf

Finally observe that the projection of $\mathcal{V}_c(e)$ to the second coordinate contains $\mathcal{T}_c(e)\setminus \mathcal{S}_c(e)$ by \cref{pr:AmazingEdge}.
Hence, the combination of \cref{pr:UpperBound,pr:LowerBound}, together with a trivial modification gives the main estimate.

\begin{theorem}\label{thm:MainEstimate}
    Let $K\in \mathbb{N}$, $e\in U_c$ be $K$-bad for $c$ and $x\in e$.
    Then we have
    $$|\fV_c(e)|\ge \frac{K^2}{(8\Delta)^4}-(16\Delta)^5K.$$
    
    In particular, for every $e\in U_c$ that is $K$-bad and $x\in e$ there are at least $\left(\frac{K^2}{(8\Delta)^4}-(16\Delta)^5K\right)$ many pairs $(f^1,f^2)$ such that $W_c(e,f^1,f^2)$ is a three step Vizing chain.
\end{theorem}

\section{Improving colorings}\label{sec:ImproveCol}

In this section we describe one step of the algorithm that will, in \cref{sec:Proof} produce the desired $\Delta(\fG)+1$ edge coloring $\mu$-almost everywhere.
The step, and therefore the whole algorithm, can be run on any Borel graph $\fG=(V,\fB,E)$ of degree bounded by $\Delta(\fG)<\infty$ endowed with an $\fE$-quasi-invariant Borel probability measure $\nu\in \fP(E)$, where $\fE=(E,\fC,I_\fG)$ is the corresponding line graph.
However, in order to show that the algorithm terminates $\nu$-almost everywhere, we need to assume additionally that $\nu$ is $\fG$-bounded, see \cref{sec:Cocycle}.

Fix $\fG$ and $\nu$ as above, and a Radon-Nikodym cocycle $\rho_\nu$ of $\nu$ with respect to $\fE$.

\begin{definition}
    We say that a partial coloring $c;E\to [\Delta(\fG)+1]$ \emph{does not admit an improvement of weight $L\in \mathbb{N}$} if 
    \begin{equation*}
        \nu\left(\left\{e\in U_c: \exists  \ 3\operatorname{-step \  Vizing \ chain} \ W_c(e) \ \operatorname{s. t.} \ \sum_{f\in W_c(e)} \rho_{\nu}(e,f)\le L\right\}\right)=0.
    \end{equation*}
    
    If this condition is not satisfied, then we say that $c$ \emph{admits an improvement of weight $L$}.
\end{definition}

\begin{theorem}\label{thm:MainStepAlgorithm}
    Let $c;E\to [\Delta(\fG)+1]$ be a partial coloring and $L\in \mathbb{N}$.
    Then there is a partial coloring $c';E\to [\Delta(\fG)+1]$ that does not admit improvement of weight $L$ with the property that $\nu(\dom(c)\setminus \dom(c'))=0$ and
    $$\nu(\{e\in E:c(e)\not=c'(e)\})\le L\nu(U_c),$$
    where $c(e)\not= c'(e)$ also includes the situation when $e\in \dom(c')\setminus \dom(c)$.
\end{theorem}
\bpf
    The strategy of the proof follows closely \cite[Proof of Proposition 5.4]{grebik2020measurable}.
    For a partial coloring $d;E\to [\Delta(\fG)+1]$ define $A_d$ to be the set of those $e\in U_d$ for which there exists a $3$-step Vizing chain $W_d(e)$ such that
    $$\sum_{f\in W_{d}(e)}\rho_\nu(e,f)\le L.$$
    Clearly, $\nu(A_d)=0$ if and only if $d$ does not admit improvement of weight $L$.
    Set $c_0=c$.
    We use induction to build a transfinite sequence of partial colorings $(c_{\alpha})_{\alpha<\aleph_1}$ that satisfy the following:
    \begin{enumerate}
        \item for every $\alpha< \beta<\aleph_1$, we have $\nu(\dom(c_\alpha)\setminus \dom(c_\beta))=0$,
        \item for every $\alpha<\aleph_1$, if $\nu(A_{c_\alpha})\not= 0$, then $\nu(U_{c_{\alpha+1}})<\nu(U_{c_{\alpha}})$,
        \item for every $\alpha<\beta<\aleph_1$, if $c_{\alpha}\not=c_{\beta}$, then $\nu(U_{c_{\beta}})<\nu(U_{c_{\alpha}})$,
        \item if $c_{\alpha}=c_{\alpha+1}$ for some $\alpha<\aleph_1$, then $c_{\alpha}=c_\beta$ for every $\alpha\le \beta<\aleph_1$,
        \item for every $\alpha< \beta<\aleph_1$,
        $$\nu(\{e\in E:c_{\alpha}(e)\not=c_\beta(e)\})\le L\sum_{\alpha\le \alpha'<\beta}\nu(S_{\alpha'})=L\nu\left(\bigcup_{\alpha\le \alpha'<\beta} S_{\alpha'}\right)\le L\nu(U_c),$$
        where $S_{\alpha'}=U_{c_{\alpha'}}\setminus U_{c_{\alpha'+1}}$.
    \end{enumerate}
    Once we build such a sequence then we are done.
    Indeed, conditions (3) and (4) guarantee the existence of $\alpha<\aleph_1$ such that $c_{\alpha}=c_{\beta}$ for every $\alpha\le \beta<\aleph_1$ as there are no strictly decreasing sequences of real numbers of length $\aleph_1$.
    Define $c'=c_{\alpha_0}$, where $\alpha_0$ is the minimal ordinal with this propery.
    Then (1) (with the choice $0< \alpha_0$) implies $\nu(\dom(c)\setminus \dom(c'))=0$, (2) implies that $c'=c_{\alpha_0}$ does not admit improvement of weight $L$ and (5) (with the choice $0< \alpha_0$) implies that $\nu(\{e\in E:c(e)\not=c'(e)\})\le L\nu(U_c)$.\\

    {\bf Successor stage $\alpha\mapsto \alpha+1$.}
    Suppose that we have constructed a sequence of partial colorings $(c_\beta)_{\beta\le \alpha}$ such that the property (1)--(5) hold for every (pair of) ordinal(s) less or equal than $\alpha$.
    If $\nu(A_{c_{\alpha}})=0$, then setting $c_{\alpha+1}=c_{\alpha}$ clearly works.
    Suppose that $\nu(A_{c_{\alpha}})>0$ and pick any Borel assignemnt $e\in A_{c_{\alpha}}\mapsto W_{c_{\alpha}}(e)$ with the property $W_{c_{\alpha}}(e)$ is a $3$-step Vizing chain and $\sum_{f\in W_c(e)}\rho_{\nu}(e,f)\le L$ for evey $e\in A_{c_{\alpha}}$.

    {\bf Case I.}
    There is $k\in \mathbb{N}$ such that 
    $$\nu(\{e\in A_{c_{\alpha}}:|W_{c_{\alpha}}(e)|=k\})>0.$$
    By \cite[Proposition~4.6]{KST} (applied on the $2k+2$ power graph of $\fE$ that is of bounded degree), there is a set $S_\alpha\subseteq A_{c_\alpha}$ with the property that
    \begin{itemize}
        \item [(a)] $\nu(S_\alpha)>0$,
        \item [(b)] $e\not=e'\in S_\alpha$ are at least $2k+2$ far apart in the graph distance of $\fE$,
        \item [(c)] $|W_{c_{\alpha}}(e)|=k$ for every $e\in S_{\alpha}$.
    \end{itemize}
    Set $T_\alpha=\bigcup_{e\in S_\alpha}W_{c_{\alpha}}(e)$ and observe that
    \begin{equation}\label{eq:Improvement1}
    \tag{*}
        \nu(T_\alpha)=\int_{e\in S_\alpha} \sum_{f\in W_{c_\alpha}(e)}\rho_{\nu} (e,f) \ d\nu \le L \nu(S_\alpha)
    \end{equation}
    by item (2) in \cref{pr:BasicCocycle}.
    As $\{W_{c_{\alpha}}(e)\}_{e\in S_\alpha}$ are pairwise of positive distance from each other by (b) and (c), and each $W_{c_{\alpha}}(e)$ is augmenting, there is a partial coloring $c_{\alpha+1}$ with the property that
    \begin{itemize}
        \item [(i)] $T_\alpha \subseteq \dom(c_{\alpha+1})$,
        \item [(ii)] $c_\alpha \upharpoonright (\dom (c_{\alpha})\setminus T_\alpha)=c_{\alpha+1}\upharpoonright (\dom (c_{\alpha})\setminus T_\alpha)$.
    \end{itemize}
    We claim that $c_{\alpha+1}$ works as required, namely, let $\alpha'<\alpha+1$, then
    \begin{enumerate}
        \item $\nu(\dom(c_{\alpha'})\setminus \dom(c_{\alpha+1}))=0$ as $\nu(\dom(c_{\alpha'})\setminus \dom(c_{\alpha}))=0$ and $\dom(c_{\alpha+1})=\dom(c_{\alpha})\cup S_{\alpha}$ by (i) and (ii),
        \item follows from $\dom(c_{\alpha+1})=\dom(c_{\alpha})\cup S_{\alpha}$ together with (a),
        \item follows from (2) combined with inductive assumption (3),
        \item if $c_{\alpha'}=c_{\alpha'+1}$ for some ${\alpha'<\alpha}$, then $\nu(A_{c_{\alpha}})=0$ by the inductive assumption (4) and (2),
        \item observe that $\nu\left(S_{\alpha}\cap S_{\beta}\right)=0$ for every $\beta<\alpha$ by the inductive assumption (1), consequently when combined with the inductive assumption (5) and \eqref{eq:Improvement1}, we have
        \begin{equation*}
            \begin{split}
                \nu(\{e\in E:c_{\alpha'}(e)\not=c_{\alpha+1}(e)\})\le & \  \nu(\{e\in E:c_{\alpha'}(e)\not=c_{\alpha}(e)\})\\
                & \ +\nu(\{e\in E:c_{\alpha}(e)\not=c_{\alpha+1}(e)\})\\
                \le & \ L\sum_{\alpha'\le \beta<\alpha}\nu(S_{\beta})+ \nu(T_{\alpha})\\
                \le & \ L\sum_{\alpha'\le \beta<\alpha+1}\nu(S_{\beta})\\
                = & \ L\nu\left(\bigcup_{\alpha'\le \beta<\alpha+1} S_{\alpha'}\right)\le L\nu(U_c).
            \end{split}
        \end{equation*}
    \end{enumerate}

    {\bf Case II.}
    There is no $k\in \mathbb{N}$ such that 
    $$\nu(\{e\in A_{c_{\alpha}}:|W_{c_{\alpha}}(e)|=k\})>0.$$
    In another words, $\nu$-almost every $3$-step Vizing chain is infinite.
    Observe that this can happen if and only if there is an assignment (defined for $\nu$-almost every $e\in A_{c_{\alpha}}$) $e\in A_{c_\alpha}\mapsto (x_e,\alpha_e,\beta_e)$, where $x_e\in V$ and $\alpha_e,\beta_e\in [\Delta+1]$, such that $W_{c_{\alpha}}(e)=M(e)^\frown P_{c_{\alpha}}(x_e,\alpha_e/\beta_e)$.
    Using finite additivity of $\nu$ and \cite[Proposition~4.6]{KST}, we find $10<k\in \mathbb{N}$, $\alpha,\beta\in [\Delta+1]$ and $R_\alpha\subseteq A_{c_\alpha}$ such that
    \begin{itemize}
        \item [(a)] $\nu(R_\alpha)>0$,
        \item [(b)] $|M(e)|=k$, $\alpha_e=\alpha$ and $\beta_e=\beta$ for every $e\in R_\alpha$,
        \item [(c)] $e\not=e'\in R_\alpha$ are at least $5k$ far apart in the graph distance of $\fE$.
    \end{itemize}
    Note that this implies that if $e\not= e'\in R_\alpha$, then
    \begin{itemize}
        \item $x_e\not=x_{e'}$ and consequently $P_{c_{\alpha}}(x_e,\alpha/\beta)$ and $P_{c_{\alpha}}(x_{e'},\alpha/\beta)$ are vertex disjoint,
        \item $M(e)$ and $M(e')$ are at least $2k$ apart in the graph distance of $\fE$.
    \end{itemize}
    However, it can happen that $M(e)\cap P_{c_{\alpha}}(x_{e'},\alpha/\beta)\not= \emptyset$.
    
    We address this issue as follows.
    Define an auxiliary directed graph $\fH$ on $R_\alpha$ as follows.
    For $e\not=e'\in R_\alpha$, let $(e,e')$ be an oriented edge if $P_{c_{\alpha}}(x_{e'},\alpha/\beta)$ intersect $B_\fE(e,2k)$, the ball of radius $2k$ around $e$.
    Note that as $|B_\fE(e,2k)|\le (2\Delta)^{2k}$ for every $e\in R_\alpha$, the graph $\fH$ has uniformly bounded outdegree.
    By \cite[Proposition~4.5]{KST}, we can write $R_\alpha=\bigcup_{n\in \mathbb{N}} R_{\alpha,n}$, where each $R_{\alpha,n}$ is $\fH$-independent.
    By $\sigma$-additivity of $\nu$ we find $n\in \mathbb{N}$ such that $\nu(R_{\alpha,n})>0$ and set $S_{\alpha}=R_{\alpha,n}$.

    Let $e\not=e'\in S_\alpha$.
    By the definition we have that $W_{c_{\alpha}}(e)$ and $W_{c_{\alpha}}(e')$ are vertex disjoint.
    Moreover, $(c_{\alpha})_e$ extends $c_{\alpha}$ as $e\in \dom((c_{\alpha})_e)$ and $\dom(c_{\alpha})\subseteq \dom((c_{\alpha})_e)$.
    Set $T_{\alpha}=\bigcup_{e\in S_\alpha} W_{c_{\alpha}}(e)$ and define
            \[
            c_{\alpha+1}(f)=
            \begin{cases}
                c_\alpha(f) & \text{if $f\not\in T_{\alpha}$} \\
                (c_\alpha)_e(f) & \text{if $f\in W_{c_{\alpha}}(e)$, where $e\in S_\alpha$ is the unique such edge}.
            \end{cases}
            \]
    It follows immediately that 
    \begin{itemize}
        \item [(i)] $T_\alpha \subseteq \dom(c_{\alpha+1})$,
        \item [(ii)] $c_\alpha \upharpoonright (\dom (c_{\alpha})\setminus T_\alpha)=c_{\alpha+1}\upharpoonright (\dom (c_{\alpha})\setminus T_\alpha)$.
    \end{itemize}
    Observe that $c_{\alpha+1}$ is a partial coloring.
    Indeed, if $x\in V$ is not covered by any edge of distance $k+2$ to some $e\in S_\alpha$, then $c_{\alpha+1}\upharpoonright N_\fG(x) =c_\alpha\upharpoonright N_\fG(x)$ (as the only modification is a shift of some of the infinite $\alpha/\beta$-paths).
    On the other hand if $x\in V$ is of distance at most $k+2$ to some $e\in S_{\alpha}$, then $c_{\alpha+1}\upharpoonright N_\fG(x)=(c_\alpha)_e\upharpoonright N_\fG(x)$, hence $c_{\alpha+1}$ is a partial coloring.
    Same reasoning as above shows that
    \begin{equation}\label{eq:Improvement}
    \tag{**}
        \nu(T_\alpha)=\int_{e\in S_\alpha} \sum_{f\in W_{c_\alpha}(e)}\rho_{\nu} (e,f) \ d\nu \le L \nu(S_\alpha)
    \end{equation}
    by item (2) in \cref{pr:BasicCocycle}.
    Verifying the conditions (1)--(5) can be done mutatis mutandis as in the Case (I).

    {\bf Limit stage $\beta\nearrow\alpha$.}
    Suppose that $\alpha<\aleph_1$ is a limit ordinal and we have constructed $(c_\beta)_{\beta<\alpha}$ such that the property (1)--(5) hold for every (pair of) ordinal(s) strictly less than $\alpha$.
    We claim that $c_{\alpha}(e):=\lim_{\beta\to\alpha} c_{\beta}(e)$ is defined $\nu$-almost everywhere, i.e., the sequence of colors $(c_\beta(e))_{\beta<\alpha}$ eventually stabilizes (in fact the colors in the sequence are changed only finitely many times) for $\nu$-almost every $e\in E$. 
    This follows from the Borel-Cantelli lemma as, by the previous paragraph, the sequence $(T_{\beta})_{\beta<\alpha}$, where $T_\beta$ is the set of edges that changed their color in the $\beta$th step, satisfies
    $$\sum_{\beta<\alpha}\nu(T_\beta)\le \sum_{\beta<\alpha}L\nu(S_\beta)\le L\nu(U_c)<\infty.$$
    Hence the set of edges $e\in E$ for which $\{\beta<\alpha:e\in T_\beta\}$ is cofinal in $\alpha$ has to be $\nu$-null.
    
    It remains to show that (1)--(5) continuous to hold with $\alpha$.
    Let $\alpha'<\alpha$, then
    \begin{enumerate}
        \item by the fact that $\alpha<\aleph_1$, the inductive assumption and construction of $c_\alpha$, we have 
        $$\nu(\{e\in \dom(c_{\alpha'}):\exists \alpha' <\beta'<\alpha \ e\not\in \dom(c_{\beta'}) \ \operatorname{or} \ \lim_{\beta\to\alpha} c_\beta(e) \ \operatorname{not \ defined}\})=0,$$
        consequently, $\nu(\dom(c_{\alpha'})\setminus \dom(c_\alpha))=0$,
        \item is not relevant in limit stages,
        \item if $c_{\alpha'}\not=c_{\alpha}$, then $c_{\alpha'}\not=c_{\alpha'+1}$ (see (4)), this implies that
        $$\nu(U_{c_{\alpha'}})>\nu(U_{c_{\alpha'+1}})\ge \nu(U_{c_\alpha})$$
        as $\nu(\dom(c_{\alpha'+1})\setminus \dom(c_\alpha))=0$ by (1) (where the strict inequality follows from the inductive assumption (3)),
        \item if $c_{\beta'}=c_{\beta''}$ for some $\beta'<\alpha$ and every $\beta'\le \beta''<\alpha$, then $c_{\alpha}(e)=\lim_{\beta\to \alpha}c_{\beta(e)}=c_{\beta'}(e)$, hence $c_{\alpha}=c_{\beta'}$,
        \item if $c_{\alpha'}(e)\not= c_\alpha(e)$, then there must be $\alpha'\le \beta<\alpha$ such that $e\in T_\beta$ as otherwise $c_\alpha(e)=\lim_{\beta\to \alpha} c_\beta(e)=c_{\alpha'}(e)$, consequently, we have
        \begin{equation*}
            \begin{split}
                \nu(\{e\in E:c_{\alpha'}(e)\not=c_\alpha(e)\})\le & \  \sum_{\alpha'\le \beta<\alpha}\nu(T_\beta)\le \sum_{\alpha'\le \beta<\alpha}L\nu(S_\beta)\\
                = & \ L\sum_{\alpha'\le \beta<\alpha}\nu(S_\beta)=L\nu\left(\bigcup_{\alpha'\le \beta<\alpha} S_\beta\right)\\
                \le  & \ L\nu\left(\bigcup_{\beta<\alpha} S_\beta\right)\le L\nu(U_c)
            \end{split}
        \end{equation*}
        by the definition of $T_\beta$ and $S_\beta$ combined with the fact that $\nu(S_\beta\cap S_{\beta'})=0$ for every $\beta<\beta'<\alpha$ which follows by the inductive assumption (1).
    \end{enumerate}
    This finishes the proof.
\epf

\section{Cocycle bounded on edges}\label{sec:Cocycle}

Recall that two Borel probablity measures $\mu,\nu\in \fP(V)$ are equivalent if $\mu(A)=0$ if and only if $\nu(A)=0$ for every $A\in\fB$.
We restate \cref{thm:CocycleIntro} in a compact form for the convenience of the reader.

\begin{theorem}\label{thm:MainCocycle}
    Let $\Delta\in \mathbb{N}$, $\fG=(V,\fB,E)$ be a Borel graph such that $\Delta(\fG)<\infty$ and $\mu\in \fP(V)$ be a $\fG$-quasi-invariant.
    Then there is an equivalent $\fG$-bounded Borel probability measure $\nu\in \mathcal{P}(V)$.
    
\end{theorem}
\bpf
    Let $\fG^{[k]}$ denote the Borel graph on $V$, where $(x,y)$ form an edge if and only if $\dist_\fG(x,y)=k$.
    Then $\fG^{[1]}=\fG$ and $\fG^{[k]}$ has degree bounded by $\Delta^k$.
    Use repeatedly \cite[Proposition~4.6]{KST} to find a Borel proper edge coloring $\{A^k_i\}_{i=1}^{2\Delta^k}$ of $\fG^{[k]}$ for each $k\in \mathbb{N}\setminus \{0\}$.
    Using any Borel linear order on $V$, vertices covered by $A^k_i$ can be split into two disjoint Borel sets $A^k_{i,0}\subseteq V$ and $A^k_{i,1}\subseteq V$ together with Borel isomorphisms $f^{k}_{i,0}:A^k_{i,0}\to A^k_{i,1}$ and $f^k_{i,1}:A^k_{i,1}\to A^k_{i,0}$ such that $\{x,f^{k}_{i,0}(x)\},\{y,f^{k}_{i,1}(y)\}\in A^{k}_{i}$ for every $x\in A^{k}_{i,0}$ and $y\in A^{k}_{i,1}$.
    We also set $A^0_{0,0}=V$ and $f^0_{0,0}=\operatorname{id}_V$.
    Observe that for every $(x,y)\in F_\fG$ there is exactly one triplet $(k,i,j)$, where $k=\dist_\fG(x,y)$, $i\in 2\Delta^{\dist_\fG(x,y)}$ and $j\in \{0,1\}$, such that $f^k_{i,j}(x)=y$.
    
    Denote as $\mu^{k}_{i,j}$ the push-forward of $\mu\upharpoonright A^k_{i,1-j}$ via $\left(f^k_{i,j}\right)^{-1}$, where $j\in \{0,1\}$.
    We have
    \begin{equation}\label{eq:DefMeasure}
        \mu^k_{i,j}(B)=\mu(f^k_{i,j}(B\cap A^k_{i,j}))=\int_{B} \one_{A^k_{i,j}}(x)\rho_\mu(x,f^k_{i,j}(x)) \ d\mu
    \end{equation}
    for every $B\in \fB$, where $\rho_\mu$ is the Radon-Nikodym cocycle with respect to $\mu$.
    In particular, $\mu^k_{i,j}(V)\le 1$.
    Define
    \begin{equation}\label{eq:GoodMeasure}
        \tilde{\nu}=\mu+\sum_{k\in \mathbb{N}\setminus \{0\}}\frac{1}{2^k} \left(\frac{1}{2\Delta^k}\sum_{i=1}^{2\Delta^k}\mu^{k}_{i,0}+\frac{1}{2\Delta^k}\sum_{i=1}^{2\Delta^k}\mu^{k}_{i,1}\right).
    \end{equation}
    \begin{claim}
        $\tilde{\nu}$ is a finite Borel measure on $V$ that is equivalent with $\mu$.
        The Radon--Nikodym derivative $\Omega=\frac{d\tilde{\nu}}{d\mu}$ can be explicitly written as
            $$\Omega(x)=1+\sum_{k\in \mathbb{N}\setminus \{0\}}\frac{1}{2^k} \left(\frac{1}{2\Delta^k}\sum_{i=1}^{2\Delta^k}\one_{A^k_{i,0}}(x)\rho_\mu(x,f^k_{i,0}(x))+\frac{1}{2\Delta^k}\sum_{i=1}^{2\Delta^k}\one_{A^k_{i,1}}(x)\rho_\mu(x,f^k_{i,1}(x))\right)$$
        for $\mu$ (and $\tilde{\nu}$) almost every $x\in V$.
        In particular, $\frac{1}{\Omega}=\frac{d\mu}{d\tilde{\nu}}$.
    \end{claim}
    \bpf
    Let $n\in \mathbb{N}$ and define 
    \begin{equation*}
        \tilde{\nu}_n=\mu+\sum_{k=1}^n\frac{1}{2^k} \left(\frac{1}{2\Delta^k}\sum_{i=1}^{2\Delta^k}\mu^{k}_{i,0}+\frac{1}{2\Delta^k}\sum_{i=1}^{2\Delta^k}\mu^{k}_{i,1}\right).
    \end{equation*}
    As $\mu^{k}_{i,j}(V)\le 1$, we see that $\tilde{\nu}_n$ is a finite Borel measure on $V$ that is equivalent with $\mu$.
    Indeed, if $\mu(A)=0$, then by the definition of $\fG$-quasi-invariance we have that $\tilde{\nu}_n(A)=0$, and if $\tilde{\nu}_n(A)=0$, then $\mu(A)=0$ as $\mu(A)\le \tilde{\nu}_n(A)$ for every $A\in \fB$.
    Moreover, it is easy to see that the Radon--Nikodym derivative $\Omega_n=\frac{d\tilde{\nu}_n}{d\mu}$ satisfies
    $$\Omega_n(x)=1+\sum_{k=1}^n\frac{1}{2^k} \left(\frac{1}{2\Delta^k}\sum_{i=1}^{2\Delta^k}\one_{A^k_{i,0}}(x)\rho_\mu(x,f^k_{i,0}(x))+\frac{1}{2\Delta^k}\sum_{i=1}^{2\Delta^k}\one_{A^k_{i,1}}(x)\rho_\mu(x,f^k_{i,1}(x))\right)$$
    by \eqref{eq:DefMeasure}.

    Observe that the limit $\Omega(x)=\lim_{n\to\infty}\Omega_n(x)$ is defined for every $x\in V$ as $\{\Omega_n(x)\}_{n\in \mathbb{N}}$ is increasing, and we have
    $$\lim_{n\to\infty}\int_{A} \Omega_n(x) \ d\mu=\int_{A} \Omega(x) \ d\mu$$
    for every $A\in \fB$ by the Monotone Convergence Theorem.
    Note that by the definition of $\Omega_n$ we have
    \begin{equation}\label{eq:RNderivative}
        \lim_{n\to\infty}\tilde{\nu}_n (A)=\int_{A} \Omega(x) \ d\mu
    \end{equation}
    for every $A\in \fB$.

    The sequence $\{\tilde\nu_n\}_{n\in \mathbb{N}}$ is a Cauchy sequence in the total variation distance $\|.\|_{TV}$.
    Indeed, we have
    $$\|\tilde{\nu}_m-\tilde{\nu}_n\|_{TV}=(\tilde{\nu}_m-\tilde{\nu}_n)(V)=\left(\sum_{k=n}^m\frac{1}{2^k} \left(\frac{1}{2\Delta^k}\sum_{i=1}^{2\Delta^k}\mu^{k}_{i,0}+\frac{1}{2\Delta^k}\sum_{i=1}^{2\Delta^k}\mu^{k}_{i,1}\right)\right)(V)\le \frac{1}{2^{n-1}}$$
    as $\tilde{\nu}_m-\tilde{\nu}_n$ is a finite Borel (positive) measure for every $m>n>1$.
    Consequently, there is a finite Borel measure $\tilde{\nu}$ on $(V,\mathcal{B})$ such that $\tilde{\nu}_n\xrightarrow{TV} \tilde{\nu}$.
    In particular, we have $\tilde{\nu}_n(A)\to \tilde{\nu}(A)<\infty$ for every $A\in\fB$.
    Combined with \eqref{eq:RNderivative}, we see that $\Omega=\frac{d\tilde{\nu}}{d\mu}$ as desired.

    It remains to show that $\mu$ and $\tilde{\nu}$ are equivalent.
    If $\mu(A)=0$, then $\tilde{\nu}(A)=\lim_{n\to\infty} \tilde{\nu}_n(A)=0$ as in that case $\tilde{\nu}_n(A)=0$ for every $n\in\mathbb{N}$.
    On the other hand, if $\tilde{\nu}(A)=0$, then $\int_A \Omega(x) \ d \mu=0$ as $\Omega=\frac{d\tilde{\nu}}{d\mu}$.
    This can only happen if $\mu(A)=0$ as $\Omega\ge 1$.
    \epf

    Let $\nu$ be a normalization of $\tilde{\nu}$, i.e., $\nu=K\tilde{\nu}$ for some $0<K<\infty$.
    It is easy to see that $\nu$ and $\tilde{\nu}$, hence also $\nu$ and $\mu$, are equivalent.
    Consequently, $\nu$ is $\fG$-quasi-invariant.
    Write $\rho_\nu$ for the Radon--Nikodym cocycle of $\nu$ with respect to $\fG$.
    \begin{claim}
    There is a $\mu$-conull $\fG$-invariant set $A\subseteq V$ such that 
    \begin{equation}
        \rho_{\nu}(x,y)=\rho_\mu(x,y)\frac{\Omega(y)}{\Omega(x)}
    \end{equation}
    for every $x,y\in A$ such that $(x,y)\in F_{\fG}$.
    \end{claim}
    \bpf
    Let $g:C\to V$ be a Borel injection such that $(x,g(x))\in F_{\fG}$ for every $x\in C$.
    Set $g(C)=D$, then $g:C\to D$ is a Borel bijection.
    As $\mu$ and $\nu$ are equivalent, $\frac{\Omega}{K}= \frac{\nu}{\mu}$ and $\rho_\mu(x,g(x))\frac{\Omega(g(x))}{\Omega(x)}$ are non-negative.
    We have
    \begin{equation*}
    \begin{split}
        \int_{C} \rho_\mu(x,g(x))\frac{\Omega(g(x))}{\Omega(x)} \ d\nu(x)= & \ \int_{C} \rho_\mu(x,g(x))\frac{\Omega(g(x))}{\Omega(x)} \frac{\Omega(x)}{K}\ d\mu(x)\\
        = & \ \int_{C} \rho_\mu(x,g(x))\frac{1}{K}\Omega(g(x)) \ d\mu(x).
    \end{split}
    \end{equation*}
    By (2) \cref{pr:BasicCocycle} applied to $f(x,y)$ that is defined as $\frac{1}{K}\Omega(y)$ whenever $g(x)=y$ and $0$ otherwise, we have
    \begin{equation*}
    \begin{split}
        \int_{C} \rho_\mu(x,g(x))\frac{1}{K}\Omega(g(x)) \ d\mu(x)= & \ \int_{D} \frac{1}{K}\Omega(y) \ d\mu(y)\\
        = & \ \nu (D)=\nu (g(C))\\
        = & \ \int_{C} \rho_{\nu}(x,g(x)) \ d\nu(x)
    \end{split}
    \end{equation*}
    and the claim follows from the uniqueness of the Radon-Nikodym cocycle.
    \epf

    It remains to show that $\rho_{\mu}(x,y)\Omega(y)\le 4\Delta \Omega(x)$ for every $(x,y)\in E$, as this clearly implies $\rho_{\nu}(x,y)=\rho_\mu(x,y)\frac{\Omega(y)}{\Omega(x)}\le 4\Delta$.
    Let $(x,y)\in E$.
    Suppose that $y\in A^k_{i,j}$, i.e., $z=f^k_{i,j}(y)$ is well-defined and satisfies $\dist_\fG(y,z)=k$.
    We already observed that the triplet $(k,i,j)$ is unique.
    As $(x,y)\in F_\fG$, there is exactly one triplet $(k',i',j')$ such that $x\in A^{k'}_{i',j'}$ and $f^{k'}_{i',j'}(x)=z$.
    Moreover, as $(x,y)\in E$ we have that $k'\le k+1$.
    Indeed, we have $\dist_\fG(x,z)=k'\le \dist_\fG(x,y)+\dist_\fG(y,z)=k+1$.
    The same reasoning with the roles of $y$ and $x$ interchanged implies that the assignment $(k,i,j)\mapsto (k',i',j')$ is a bijection.
    Observe that
    \begin{equation}\label{eq:estimate}
        \frac{1}{2^{k+1}(2\Delta^{k+1})}\one_{A^k_{i,j}}(y)\rho_\mu(x,f^k_{i,j}(y))\le \frac{1}{2^{k'}(2\Delta^{k'})}\one_{A^{k'}_{i',j'}}(x)\rho_\mu(x,f^{k'}_{i',j'}(x))
    \end{equation}
    holds whenever $y\in A^k_{i,j}$.
    We have
    \begin{equation*}
    \begin{split}
        \rho_\mu(x,y)\Omega(y)= & \ \rho_\mu(x,y)\Bigg(1+\sum_{k\in \mathbb{N}\setminus \{0\}}\frac{1}{2^k} \Bigg(\frac{1}{2\Delta^k}\sum_{i=1}^{2\Delta^k}\one_{A^k_{i,0}}(y)\rho_\mu(y,f^k_{i,0}(y))\\
        & \ +\frac{1}{2\Delta^k}\sum_{i=1}^{2\Delta^k}\one_{A^k_{i,1}}(y)\rho_\mu(y,f^k_{i,1}(y))\Bigg)\Bigg)\\
        = & \ \rho_\mu(x,y)+\sum_{k\in \mathbb{N}\setminus \{0\}}\frac{1}{2^k} \Bigg(\frac{1}{2\Delta^k}\sum_{i=1}^{2\Delta^k}\one_{A^k_{i,0}}(y)\rho_\mu(x,f^k_{i,0}(y))\\
        & \ +\frac{1}{2\Delta^k}\sum_{i=1}^{2\Delta^k}\one_{A^k_{i,1}}(y)\rho_\mu(x,f^k_{i,1}(y))\Bigg)\\
        \le & \ 4\Delta \Bigg(\frac{1}{4\Delta}\rho_\mu(x,f^0_{0,0}(y))+ \sum_{k\in \mathbb{N}\setminus \{0\}}\frac{1}{2^{k+1}} \Bigg(\frac{1}{2\Delta^{k+1}}\sum_{i=1}^{2\Delta^k}\one_{A^k_{i,0}}(y)\rho_\mu(x,f^k_{i,0}(y))\\
        & \  +  \frac{1}{2\Delta^{k+1}}\sum_{i=1}^{2\Delta^k}\one_{A^k_{i,1}}(y)\rho_\mu(x,f^k_{i,1}(y))\Bigg)\Bigg)\\   
        \le & \ 4\Delta\Bigg(1+\sum_{k'\in \mathbb{N}\setminus \{0\}}\frac{1}{2^{k'}} \Bigg(\frac{1}{2\Delta^{k'}}\sum_{i=1}^{2\Delta^{k'}}\one_{A^{k'}_{i',0}}(x)\rho_\mu(x,f^{k'}_{i',0}(x))\\
        & \ +\frac{1}{2\Delta^{k'}}\sum_{i=1}^{2\Delta^{k'}}\one_{A^{k'}_{i',1}}(x)\rho_\mu(x,f^{k'}_{i',1}(x))\Bigg)\Bigg)\\
        = & \ 4\Delta \Omega(x),
    \end{split}
    \end{equation*}
    where we used (1) \cref{pr:BasicCocycle} to get the second equality and \eqref{eq:estimate} together with the fact that $(k,i,j)\mapsto (k',i',j')$ is a bijection to get the second inequality.
\epf

\section{Double counting argument} \label{sec:DoubleCount}

In this section we show that if a partial edge coloring $c$ does not admit improvement of weight $L$, then the measure of uncolored edges has to be $O\left(\frac{1}{\log^2 (L)L}\right)$ under the assumption that $\nu$ is $\fG$-bounded.

\begin{theorem}\label{thm:MainExpansion}
    Let $\fG=(V,\fB,E)$ be a Borel graph such that $\Delta(\fG)<\infty$, $\fE=(E,\fC,I_\fG)$ be the corresponding line graph, $\nu\in \fP(E)$ be $\fE$-bounded, i.e., $\nu$ satisfy $\rho_\nu(e,f)\le 8\Delta$ for every $e,f\in E$ such that $(e,f)\in I_\fG$, and $c;E\to [\Delta(\fG)+1]$ be a partial coloring that does not admit improvement of weight $L\in \mathbb{N}$, where $\log_{8\Delta}(L)\ge (8\Delta)^{20}$.
    Then
    $$\nu(U_c)\le \frac{64 (4\Delta)^7(\Delta!)^{14}}{\log_{8\Delta}^2(L) L}.$$
\end{theorem}

\begin{proposition}
    Let $\fG=(V,\fB,E)$ be a Borel graph such that $\Delta(\fG)<\infty$, $\fE=(E,\fC,I_\fG)$ be the corresponding line graph and $\nu\in \fP(E)$ be $\fE$-bounded.
    Suppose that $c;E\to [\Delta(\fG)+1]$ is a partial coloring that does not admit improvement of weight $L\in \mathbb{N}$, where $\frac{\log_{8\Delta}(L)}{2}-1\ge 2\Delta$.
    Then $\nu$-almost every $e\in U_c$ is $\frac{\log_{8\Delta}(L)}{4}$-bad for $c$.
\end{proposition}
\bpf
    Let $W_c(e)$ be a $3$-step Vizing chain at $e$.
    We have
    $$L\le \sum_{f\in W_c(e)} \rho_\nu(e,f)\le \sum_{k\le |W_c(e)|}(8\Delta)^k\le (8\Delta)^{|W_c(e)|+1}.$$
    Consequently, $|W_c(e)|\ge \log_{8\Delta}(L)-1\ge \frac{\log_{8\Delta}(L)}{2}+2\Delta$.
\epf

We get an immediate corollary of \cref{thm:MainEstimate}.

\begin{proposition}\label{pr:CorollaryEstimate}
    Let $\fG=(V,\fB,E)$ be a Borel graph such that $\Delta(\fG)<\infty$, $\fE=(E,\fC,I_\fG)$ be the corresponding line graph and $\nu\in \fP(E)$ be $\fE$-bounded.
    Suppose that $c;E\to [\Delta(\fG)+1]$ is a partial coloring that does not admit improvement of weight $L\in \mathbb{N}$, where $\log_{8\Delta}(L)\ge (8\Delta)^{20}$.
    Then we have 
    $$|\fV_c(e)|\ge \frac{1}{(4\Delta)^7}\log^2_{8\Delta}(L)$$
    for $\nu$-almost every $e\in U_c$.
\end{proposition}

\begin{proposition}\label{pr:EstimateOnThirdPath}
     Let $\fG=(V,\fB,E)$ be a Borel graph such that $\Delta(\fG)<\infty$, $\fE=(E,\fC,I_\fG)$ be the corresponding line graph and $\nu\in \fP(E)$ be $\fE$-bounded.
    Suppose that $c;E\to [\Delta(\fG)+1]$ is a partial coloring that does not admit improvement of weight $L\in \mathbb{N}$, where $\log_{8\Delta}(L)\ge (8\Delta)^{20}$.
    Then for $\nu$-almost every $e\in U_c$ and every $(f^1,f^2)\in \fV_c(e)$ we have
    $$\sum_{f\in P^3_c}\rho_\nu(e,f)\ge \frac{L}{2},$$
    where $W_c(e,f^1,f^2)=\left(F^1_c\right)^\frown \left(P^1_c\right)^\frown\left(F^2_c\right)^\frown\left(P^2_c\right)^\frown\left(F^3_c\right)^\frown\left(P^3_c\right)$.
\end{proposition}
\bpf
    Set $P=\left(F^1_c\right)^\frown \left(P^1_c\right)^\frown\left(F^2_c\right)^\frown\left(P^2_c\right)^\frown\left(F^3_c\right)$.
    It follows from the definition of $\fV_c(e)$ that $s=l(P)\le 3\Delta+\frac{\log_{8\Delta}(L)}{2}$.
    We have
    \begin{equation*}
        \begin{split}
            L\le \sum_{f\in W_c(e,f^1,f^2)}\rho_\nu(e,f)= & \  \sum_{f\in P}\rho_\nu(e,f)+\sum_{f\in P^3_c}\rho_\nu(e,f)\\
            \le & \ (8\Delta)^{3\Delta+\frac{\log_{8\Delta}(L)}{2}+2} + \sum_{f\in P^3_c}\rho_\nu(e,f).
        \end{split}
    \end{equation*}
    Consequently, $\frac{L}{2}\le L-(8\Delta)^{3\Delta+2} L^{\frac{1}{2}}\le \sum_{f\in P^3_c}\rho_\nu(e,f)$ as needed.
\epf

Define an auxiliary Borel oriented bipartite multi-graph $\fH_c$ with vertex set $U_c\sqcup \dom(c)$ such that $(e,f)$ is an edge if $f\in P^3_c$ for some
$$W_c(e,f^1,f^2)=\left(F^1_c\right)^\frown \left(P^1_c\right)^\frown\left(F^2_c\right)^\frown\left(P^2_c\right)^\frown\left(F^3_c\right)^\frown\left(P^3_c\right),$$
where $(f^1,f^2)\in \fV_c(e)$.
Note that $\fH_c$ is in general a multi-graph as there might be different $3$-step Vizing chains at the same edge $e$ for which $(e,f)\in \fH_c$ for the same $f\in \dom(c)$.
The following is the crucial observation for the double counting argument, note that the right-hand side of the inequality does not depend on $L$.

\begin{proposition}\label{cl:DoubleCounting}
    $\deg_{\fH_c}(f)\le 32 (\Delta!)^{14}$ for every $f\in \dom(c)$.
\end{proposition}
\bpf
    There are at most $2\Delta$ choices for $\alpha_3,\beta_3$ and $z_3$ such that $f\in P^3_c=P_c(z_3,\alpha_3/\beta_3)$.
    There are at most $\Delta$ choices for $y_3$, and at most $\Delta!\Delta$ choices for a fan $F^3_c$ at $x_3$.
    For $P^2_c$ we deduce that there are at most $4(\Delta+1)^2$ choices for $z_2$ and $\alpha_2,\beta_2$ as the last edge of $P^2_c$ has to intersect the first edge of $F^3_c$.
    Similar estimates hold for $F^2_c$, $P^1_c$ and $F^1_c$.
    Altogether we get that
    $$\deg_{\fH_c}(f)\le 2\Delta(\Delta! \Delta^2)^3(4(\Delta+1)^2)^2\le 32 (\Delta!)^{14}$$
    as desired.
\epf

\bpf[Proof of \cref{thm:MainExpansion}]
    Define a function ${\bf F}(e,f)$ that counts the number of oriented edges from $e$ to $f$ in $\fH_c$.
    By (2)~\cref{pr:BasicCocycle}, we have
    \begin{equation}\label{eq:doublecounting}
        \tag{DC}
        \int_{E}\sum_{f\in [e]_\fE}{\bf F}(e,f)\rho_{\nu}(e,f) \ d\nu(e)=\int_{E}\sum_{e\in [f]_\fE}{\bf F}(e,f) \ d\nu(f).
    \end{equation}
    Using \cref{cl:DoubleCounting}, we get an upper bound for the right-hand side of \eqref{eq:doublecounting} as
    $$\int_{E}\sum_{e\in [f]_\fE}{\bf F}(e,f) \ d\nu(f)\le \int_{\dom(c)}\deg_{\fH_c}(f) \ d\nu(f)\le 32 (\Delta!)^{14}.$$
    Using the definition of $\fV_c(e)$ for $e\in U_c$, \cref{pr:EstimateOnThirdPath} and \cref{pr:CorollaryEstimate}, we get an lower bound for the left-hand side of \eqref{eq:doublecounting} as
    \begin{equation*}
        \begin{split}
            \int_{E}\sum_{f\in [e]_\fE}{\bf F}(e,f)\rho_{\nu}(e,f) \ d\nu(e)= & \ \int_{U_c}\sum_{(f^1,f^2)\in \fV_c(e)}\sum_{f\in P^3_c}\rho_{\nu}(e,f) \ d\nu(e)\\
            \ge & \ \int_{U_c}\sum_{(f^1,f^2)\in \fV_c(e)}\frac{L}{2} \ d\nu(e)\\
            \ge & \ \int_{U_c} \frac{1}{(4\Delta)^7}\log^2_{8\Delta}(L)\frac{L}{2} \ d\nu(e)=\frac{1}{2(4\Delta)^7}\log^2_{8\Delta}(L)L\nu(U_c).
        \end{split}
    \end{equation*}
    Altogether, we infer that
    $$\nu(U_c)\le \frac{64 (4\Delta)^7(\Delta!)^{14}}{\log^2_{8\Delta}(L)L}$$
    as desired.
\epf

\section{Proof of the main result}\label{sec:Proof}

We restate \cref{thm:MainIntro}, in a compact form, for the convenience of the reader.

\begin{theorem}
    Let $\fG=(V,\fB,E)$ be a Borel graph such that $\Delta(\fG)<\infty$ and $\mu\in \fP(V)$.
    Then there is a Borel map $c:E\to [\Delta(\fG)+1]$ that is a proper edge coloring $\mu$-almost everywhere.
\end{theorem}
\bpf
    Apply \cref{pr:BasicMeasure} to $\fG$ to get $\hat\mu$, then apply \cref{thm:MainCocycle} to $\fE=(E,\fC,I_\fG)$ and $\hat\mu$ to get $\nu\in \fP(E)$ with the property that $\rho_\nu(e,f)\le 8\Delta$ for every $e,f\in E$ such that $(e,f)\in I_\fG$, and $\mu(\{x\in V:\exists e\in A \ x\in e\})=0$ for every $A\in \fC$ such that $\nu(A)=0$.
    
    Set $L_n=A\cdot (8\Delta)^n$, where $A\in \mathbb{N}$ is such that $\log_{8\Delta}(L_n)\ge (8\Delta)^{20}$ for every $n\in \mathbb{N}$, and define inductively, using \cref{thm:MainStepAlgorithm}, a sequence $\{c_n\}_{n\in \mathbb{N}}$ such that $c_n$ does not admit improvement of weight $L_n$ for every $n\in \mathbb{N}$ and
    $$\nu(\{e\in E:c_n(e)\not=c_{n+1}(e)\})\le L_{n+1}\nu(U_c)$$
    for every $n\in \mathbb{N}$.

    The Borel--Cantelli Lemma implies that $c(e)=\lim_{n\to\infty} c_n(e)$ is defined $\nu$-almost everywhere as
    \begin{equation*}
        \begin{split}
            \sum_{n\in \mathbb{N}}L_{n+1}\nu(U_{c_n})\le & \   \sum_{n\in \mathbb{N}}\frac{64 (4\Delta)^7(\Delta!)^{14}L_{n+1}}{\log_{8\Delta}^2(L_n)L_n}\\
            = & \ \sum_{n\in \mathbb{N}}\frac{64 (4\Delta)^7(\Delta!)^{14}(8\Delta)^{n+1}}{n^2 (8\Delta)^{n}}=\sum_{n\in \mathbb{N}}\frac{128 (4\Delta)^8(\Delta!)^{14}}{n^2}<\infty
        \end{split}
    \end{equation*}
    by \cref{thm:MainExpansion}.
    Altogether, this shows that $c$ is correct at $\mu$-almost every vertex by the definition of $\nu$.    
\epf

\bibliographystyle{alpha}
\bibliography{ref}

\end{document}